\documentclass[12pt]{article}

\title{Galois Extensions of Braided Tensor Categories \\ and Braided Crossed G-Categories}

\author{Michael M\"uger\thanks{Supported by NWO.} \\
Faculteit Wiskunde en Informatica, Universiteit Utrecht, Netherlands \\
and Korteweg-de Vries Institute, Amsterdam, Netherlands \\ 
email: {\tt mmueger@science.uva.nl}}

\newlength{\dinwidth}
\newlength{\dinmargin}
\usepackage{latexsym, amssymb, amsmath, theorem,tan-gle}
\usepackage{tan-gle}

\input diagrams         

\def\1#1{{\bf #1}}
\def\2#1{{\cal #1}}
\def\3#1{{\sl #1}}
\def\4#1{{\tt #1}}
\def\5#1{{\sf #1}}
\def\6#1{{\mathfrak #1}}
\def\7#1{{\mathbb #1}}

\newcommand{\DS}{\displaystyle}
\newcommand{\ba}{\begin{array}}
\newcommand{\ea}{\end{array}}
\newcommand{\bea}{\begin{eqnarray}}
\newcommand{\eea}{\end{eqnarray}}
\newcommand{\bean}{\begin{eqnarray*}}
\newcommand{\eean}{\end{eqnarray*}}
\newcommand{\nn}{\nonumber}

\newcommand{\ve}{\varepsilon}

\newcommand{\impl}{\Rightarrow}
\newcommand{\rarr}{\rightarrow}

\newcommand{\ol}{\overline}

\newcommand{\del}{\partial}
\newcommand{\id}{\mathrm{id}}
\newcommand{\mcirc}{\,\circ\,}

\newcommand{\Hom}{\mathrm{Hom}}
\newcommand{\End}{\mathrm{End}}
\newcommand{\Aut}{\mathrm{Aut}}

\newcommand{\Mod}{\mathrm{Mod}}
\newcommand{\Obj}{\mathrm{Obj}\,}
\newcommand{\Rep}{\mathrm{Rep}}
\newcommand{\GLoc}{G\!-\!\mathrm{Loc}}

\def\endexem{\hfill{$\Box$}\medskip}
\newcommand{\qed}{\ \hfill $\blacksquare$\medskip}

\theoremstyle{change}
\theoremheaderfont{\scshape}
\newtheorem{defin}{Definition}[section]
\newtheorem{defprop}[defin]{Definition/Proposition}
\newtheorem{lemma}[defin]{Lemma}
\newtheorem{prop}[defin]{Proposition}
\newtheorem{theorem}[defin]{Theorem}
\newtheorem{coro}[defin]{Corollary}
\newtheorem{conj}[defin]{Conjecture}
\theorembodyfont{\rmfamily}
\newtheorem{rema}[defin]{Remark}
\newtheorem{noname}[defin]{}

\newcommand{\bdefin}{\begin{defin}}
\newcommand{\bdefprop}{\begin{defprop}}
\newcommand{\blemma}{\begin{lemma}}
\newcommand{\bprop}{\begin{prop}}
\newcommand{\btheor}{\begin{theorem}}
\newcommand{\bcoro}{\begin{coro}}
\newcommand{\edefin}{\end{defin}}
\newcommand{\edefprop}{\end{defprop}}
\newcommand{\elemma}{\end{lemma}}
\newcommand{\eprop}{\end{prop}}
\newcommand{\etheor}{\end{theorem}}
\newcommand{\ecoro}{\end{coro}}
\newcommand{\bconj}{\begin{conj}}
\newcommand{\econj}{\end{conj}}
\newcommand{\brem}{\begin{rema}}
\newcommand{\erem}{\endexem\end{rema}}
\newcommand{\bnix}{\begin{noname}}
\newcommand{\enix}{\end{noname}}

\newcommand{\prf}{\noindent{\it Proof. }}

\def\mobj#1{\raise .4\unitlens\hbox{\put(0,0){$#1$}}}

\def\mychi{\raise 2pt\hbox{$\chi$}}

\begin{document}
\maketitle\noindent

\numberwithin{equation}{section}

\abstract{We show that the author's notion of Galois extensions of braided tensor
categories \cite{mue06}, see also \cite{brug1}, gives rise to braided crossed
G-categories, recently introduced for the purposes of 3-manifold topology \cite{t2}. 
The Galois extensions $\2C\rtimes\2S$ are studied in detail, in particular we determine 
for which $g\in G$ non-trivial objects of grade $g$ exist in $\2C\rtimes\2S$.
}


\section{Introduction}
According to the influential paper \cite{js}, the notion of braided tensor categories
(btc for short) originated in (I) considerations in higher dimensional category theory
(btc as 3-categories with one object and one 1-morphism) and (II) homotopy theory (braided
categorical groups classifying connected homotopy types with only $\pi_2, \pi_3$
non-trivial). On the other hand, the (III) representation categories of quasitriangular
(quasi-, weak etc.) Hopf algebras, cf.\ e.g.\ \cite{ka}, and of (IV) quantum field
theories (QFT) in low-dimensional space times \cite{frs1,fk}, in particular conformal
field theories \cite{frs1,ms}, are btc. Finally, (V) the category of tangles is a btc, which is
the origin of various constructions of invariants of links and 3-manifolds
\cite{t1,ka,baki}. It goes without saying that all five areas continue to be very active
fields of research and the connections continue to be explored. 

In this paper we are concerned with a recent generalization of the notion of btc which is
quite interesting in that can be approached from most of the above viewpoints. (V): In the
context of his programme of homotopy TQFT, Turaev \cite{t1} introduced {\bf braided
G-crossed categories} and showed that, subject to some further conditions, they give
rise to invariants of 3-dimensional G-manifolds, to wit 3-manifolds together with a
principal G-bundle. Let us state the definition in its simplest form.

\bdefin \label{def-crossedG}
Let $G$ be a (discrete) group. A strict crossed G-category is a strict tensor category
$\2C$ together with 
\begin{itemize}
\item a map $\del: \Obj\,\2C\rarr G$ constant on isomorphism classes,
\item a homomorphism $\gamma: G\rarr\Aut\,\2C$ (strict monoidal automorphisms of $\2C$)
\end{itemize}
such that 
\begin{enumerate}
\item $\del(X\otimes Y)=\del X\,\del Y$.
\item $\del(\gamma_g(X))=g\,\del(X)\,g^{-1}$.
\end{enumerate}
We write ${}^Y\!\!\cdot=\gamma_{\del Y}(\cdot)$. A braiding for a crossed G-category $\2C$ is a
family of isomorphisms $c_{X,Y}: X\otimes Y\rarr{}^X\!Y\otimes X$ such that 
\bean c_{X,Z\otimes T} &=& \id_{{}^X\!Z}\otimes c_{X,T}\mcirc
          c_{X,Z}\otimes\id_T, \\
   c_{X\otimes Y,Z} &=& c_{X,{}^Y\!Z}\otimes\id_Y\mcirc\id_X\otimes c_{Y,Z}, 
\eean
\[ c_{X',Y'}\mcirc s\otimes t={}^X\!t\otimes s\mcirc c_{X,Y} \quad\forall
   s: X\rarr X', \ t: Y\rarr Y'. \]
\edefin
Of the various generalizations permitted by this definition we will need only the
admission of inhomogeneous objects, cf.\ Section \ref{sec3}. As to subject (III): In
\cite{t3} it was shown that some crossed G-categories can be obtained from quantum
groups. With a view towards applications to algebraic topology (II), in \cite{carr} a
notion of {\bf categorical G-crossed module} was defined. The latter are simply crossed
G-categories which are categorical groups, i.e.\ monoidal groupoids with invertible
objects. In turn, categorical G-crossed modules generalize Whitehead's ubiquitous notion
of crossed modules and Conduch\'e's 2-crossed modules. 

The main result of the present paper is to show that braided crossed G-categories arise 
from a categorical construction, the Galois extensions of braided tensor categories 
\cite{mue06,brug1}. This refers to the construction in \cite{mue06}
which associates to a braided tensor category $\2C$ and a full symmetric subcategory $\2S$
a tensor category $\2C\rtimes\2S$. The braiding of $\2C$ lifts to a braiding of
$\2C\rtimes\2S$ iff $\2S$ is contained in the center $Z_2(\2C)$ of $\2C$, the latter being
the full subcategory of objects $X$ satisfying $c_{X,Y}\circ c_{Y,X}=\id$ for all $Y\in\2C$.
(In \cite{brug1}, where a category equivalent to $\2C\rtimes\2S$ was defined, the objects
of $Z_2(\2C)$ were called transparent.) Dropping the condition $\2S\subset Z_2(\2C)$ we
show in Theorem \ref{t-main1} that $\2C\rtimes\2S$ is a braided crossed G-category, where
we also clarify for which $g\in G$ there exist $X_g\in\2C\rtimes\2S$ with $\del X=g$, cf.\
Theorem \ref{theor-spec}. In the final Subsection \ref{ss-34} we show that a subcategory
$\2S\subset\2C$ where $\2S\cong\Rep\,G$ with $G$ finite abelian induces a $G$-grading on
$\2C$ compatible with the one on $\2C\rtimes\2S$. Similar results are obtained in
\cite{ki}, in particular part II. However, our approach is quite different, more suitable
for the application to quantum field theory \cite{mue15} sketched below, and in places
somewhat more satisfactory, e.g.\ concerning the braiding on $\2C\rtimes\2S$.

We close this introduction with a glance at the applications of this paper in quantum
field theory and topology.
In a companion paper \cite{mue15} we will show, in the context of algebraic quantum field
theory \cite{haag}, that a chiral conformal field theory $A$ carrying an action of a
finite group $G$ gives rise to a braided crossed G-category $\GLoc\,A$ of `G-twisted 
representations'. The full subcategory $\del^{-1}(e)\subset\GLoc\,A$ of grade zero objects
is just the ordinary braided representation category $\Rep\,A$, which does
not the G-action into account. In \cite{mue15} we prove the equivalences
\bean \GLoc\,A &\simeq& \Rep\,A^G\rtimes\2S, \\
   \Rep\,A^G &\simeq& (\GLoc\,A)^G, \eean
where $A^G$ is the `orbifold theory' \cite{dvvv}, i.e.\ the subtheory of $A$ consisting of the
fixpoints under the G-action, and $\2S\simeq\Rep\,G$ is a full subcategory of
$\Rep\,A^G$. The significance of the first equivalence is that the same braided crossed
G-category arises (i) as the -- intrinsically defined -- category of G-twisted
representations of $A$ and (ii) by the crossed product construction of \cite{mue06} whose
braided crossed G-structure is the subject of the present work. The second equivalence
computes the representation category of the orbifold theory $A^G$ in terms of $\GLoc\,A$,
i.e.\ categorical information about $A$. To put this into context we emphasize the well
known fact that the grade zero subcategory $\Rep\,A\subset \GLoc\,A$ does {\it not}
contain enough information to determine $\Rep\,A^G$.

Finally, by \cite{klm} the categories $\Rep\,A$ and $\Rep\,A^G$ are modular, and Corollary
\ref{coro-spec} implies that $\GLoc\,A\simeq\Rep\,A^G\rtimes\2S$ has full G-spectrum, i.e.\
there exists an object of grading $g$ for every $g\in G$. Combining this with Turaev's
work \cite{t1, t2} on invariants of (G-)manifolds we thus obtain an equivariant version of
the chain 
\[ \mathrm{Rational\ chiral\ CFT}\ \leadsto\ \mathrm{modular\ category}\ \leadsto\
  \mathrm{3-manifold\ invariant}, \]
of constructions, namely
\[ \parbox{3.5cm}{Rational chiral CFT with symmetry G}\ \ \leadsto \ \
  \parbox{2.7cm}{modular crossed G-category}\ \ \leadsto\ \ 
  \parbox{6cm}{invariant for 3-manifolds equipped with principal G-bundle}. \]

The above applications of the constructions of this paper place braided crossed
G-categories squarely into the context of the areas (I) (higher category theory) and (V) 
(quantum field theory) mentioned above. Most results of this paper and of \cite{mue15}
were announced in \cite{mue08}.


\section{Preliminaries}
We briefly recall without proof the facts concerning Tannakian and module categories that
will be needed later. Some of those are well known, while others are relatively recent.

We assume as known the notions of abelian, monoidal (or tensor) braided, symmetric, rigid
and ribbon categories, cf.\ e.g.\ \cite{cwm,js,ka,baki}. All categories considered in
this paper will be $\7F$-linear semisimple (thus in particular abelian) over an
algebraically closed field $\7F$ with finite dimensional $\Hom$-spaces and monoidal with
$\End\,\11=\7F\id_\11$. Unless otherwise stated tensor categories will be 
strict, as we are allowed to assume by virtue of the coherence theorems. A $\7C$-linear
tensor category is a $*$-category if there exists a $*$-operation, i.e.\ an involutive
antilinear contravariant and monoidal endofunctor $*$ that acts trivially on the
objects. In other words, $s^*\in\Hom(Y,X)$ if $s\in\Hom(X,Y)$, $s^{**}=s$ and, whenever
these expressions are defined, $(s\circ t)^*=t^*\circ s^*$ and 
$(s\otimes t)^*=s^*\otimes t^*$. A $*$-operation is positive if $s^*\circ s=0$ implies
$s=0$. A category with positive $*$-operation is called $*$-category \cite{glr,dr6,lro} or 
unitary \cite{t1}, cf.\ also \cite{y3}. (Since we assume finite dimensional hom-spaces, a
$*$-category in fact is a $C^*$- and $W^*$-category in the sense of \cite{glr,dr6},
cf.\ e.g.\ \cite[Proposition 2.1]{mue06}.)

The category of finite dimensional polynomial representations of a reductive proalgebraic
group (in characteristic zero) is a rigid abelian symmetric tensor category with
$\End\,\11=\7F\id_\11$. The category of finite dimensional continuous representations of a
compact topological group has the same properties and is in addition a $*$-category. There
are converses to these statements due to Doplicher and Roberts \cite{dr6} and to Deligne
\cite{del}, respectively. For our purposes in this paper it is sufficient to consider
symmetric categories with finitely many (isomorphism classes of) simple objects,
corresponding to finite groups. 

\bdefin
\begin{enumerate}
\item A TC is a semisimple $\7F$-linear spherical tensor category \cite{bw} with finite
dimensional $\Hom$-spaces and $\End\,\11=\7F\id_\11$, where $\7F$ is an algebraically closed
field. It is called finite if the set of isomorphism classes of simple objects is finite.
The dimension of a finite TC is given by $\dim\2C=\sum_i d(X_i)^2$, where $i$ runs
through the set of isomorphism classes of simple objects and $d$ is the dimension function
defined by the spherical structure.
\item A BTC is a semisimple $\7F$-linear rigid braided ribbon category with finite
dimensional $\Hom$-spaces and $\End\,\11=\7F\id_\11$, and is automatically a TC. 
\item An STC is a symmetric BTC.
\item An STC over $\7F$ is admissible if either (i) $\7F=\7C$, $\2C$ is a $*$-category and all
objects have trivial twist $\Theta(X)$, or (ii) $\7F$ has characteristic zero and
$d(X)\in\7Z_+$ for all $X\in\2C$.  
\end{enumerate}
\edefin

\brem Ad 1: Since we work over algebraically closed fields throughout, an object $X$ is
simple (every non-zero subobject is isomorphic to $X$) iff it is absolutely simple
($\End\,X=\7F\id_X$). We will therefore just speak of simple objects.

By dropping the assumption of sphericity one arrives at the notion of fusion
categories which were studied in \cite{eno}. There are remarkably strong results like the
automatic positivity of $\dim\2C$ when $\7F=\7C$. (Yamagami has shown \cite{y3} that a
$*$-structure gives rise to an essentially unique spherical structure, and one might
suspect that this generalizes to fusion categories.)

Ad 2: A rigid ribbon category gives rise to a spherical structure and conversely in a
spherical braided category $\2C$ there exists a canonical twist $\Theta$ rendering $\2C$ a
ribbon category. See \cite{bw,ye}.

Ad 3: At first sight, the supplementary conditions (i) and (ii) on the twists and the
dimensions, respectively, look quite different. This is due to the different notions of
duality in both formalisms, but ultimately both conditions amount to the same thing.
Let $X\in\2C$. In \cite{dr6} one chooses 
$r_X: \11\rarr \ol{X}\otimes X, \ol{r}_X: \11\rarr X\otimes \ol{X}$ such that 
\bean \id_{\ol{X}}\otimes \ol{r}_X^* \mcirc r_X\otimes\id_{\ol{X}} &=& \id_{\ol{X}}, \\
  \id_X\otimes r_X^* \mcirc\ol{r}_X\otimes\id_X &=& \id_X, \\
  r_X^*\circ r_X = \ol{r}_X^*\circ\ol{r}_X &=& d(X)\id_\11. 
\eean
One then defines the twist $\Theta(X)\in\End\,X$ by 
\[ \Theta(X)=r_X^*\otimes\id_X\mcirc\id_{\ol{X}}\otimes c_{X,X}\mcirc r_X\otimes\id_X.\]
For simple $X$ one finds $\Theta(X)=\pm\id_X$, whereas $d(X)\ge 0$ is automatic by
positivity of the $*$-operation. In fact, one proves $d(X)\in\7Z_+$, and the condition
$\Theta(X)=\id_X$ is necessary and sufficient for $\2C\simeq\Rep\,G$ for some $G$.

On the other hand, in \cite{del} one has morphisms $d_X: \11\rarr X\otimes\ol{X}$, 
$e_X: \ol{X}\otimes X\rarr\11$, which are part of the given data and satisfy the usual
triangular equations. One then defines
\[ \delta_X=c_{X,\ol{X}}\circ d_X: \11\rarr\ol{X}\otimes X, \quad
   \eta_X=e_X\circ c_{X,\ol{X}}: X\otimes\ol{X} \rarr\11. \]
With this definition the twist 
$\Theta(X)=e_X\otimes\id_X\mcirc\id_{\ol{X}}\otimes c_{X,X}\mcirc \delta_X$ is
automatically trivial, but $d(X)=\eta_X\circ d_X=e_X\circ\delta_X$ is not necessarily
positive. In any case, for a $*$-category one has both notions of duals, and the
supplementary conditions are equivalent.
\erem

\btheor \label{theor2} \cite{dr6,del} 
Let $\2C$ be a finite admissible STC over $\7F$. Then there exists a finite group $G$,
unique up to isomorphism, such that there is an equivalence $\2C\simeq \Rep_\7F^{fin}\,G$
compatible with all structures in sight. 
\etheor

\brem 1. The proof of Theorem \ref{theor2} in \cite{del} roughly consists of two steps:
(i) One constructs a faithful tensor functor $E: \2C\rarr\mathrm{Vect}_\7F^{fin}$. (ii)
Defining $G=\mathrm{Nat}^\otimes E$, the set of monoidal natural transformations from $E$
to itself, one finds (a) $G$ is a group, cf.\ e.g.\ \cite[Proposition 7.1]{js}, by virtue
of rigidity of $\2C$, and (b) $\2C\simeq\Rep_\7F\,G$. 

2. If $\2C$ has objects with non-trivial twists or integral but non-positive dimensions,
respectively, it still is the representation category of a supergroup, i.e.\ a pair
$(G,k)$ where $G$ is a group and $k\in Z(G)$ is involutive, cf.\ \cite[Section 7]{dr6},
see also \cite{del2}. This generalization will not be used in this paper.
\erem

\bdefin \label{def_frob}
Let $\2C$ be a strict tensor category. A Frobenius algebra in $\2C$ is a quintuple
$(\Gamma,m,\eta,\Delta,\ve)$ such that 
$(\Gamma\in\2C, m:\Gamma^2\rarr\Gamma,\eta:\11\rarr\Gamma)$ is a monoid,
$(\Gamma,\Delta: \Gamma\rarr\Gamma^2, \ve:\Gamma\rarr\11)$ is a comonoid and the condition
\[ \id_\Gamma\otimes m\mcirc\Delta\otimes\id_\Gamma \ = \ \Delta\circ m \ =\
  \Delta\otimes\id_\Gamma\mcirc \id_\Gamma\otimes m \]
holds.
A Frobenius algebra in an $\7F$-linear category is called strongly separable \cite{mue09} if
\[ m\circ\Delta=\alpha\,\id_\Gamma, \quad \ve\circ\eta=\beta\,\id_\11, \quad \alpha,\beta\in
   \7F^*. \]
\edefin

\brem Following earlier terminology used by the author, which in turn was inspired by 
F.\ Quinn, strongly separable Frobenius algebras were called `special' in \cite{fs}.
\erem

\bprop \label{prop2}\cite{mue09}
Let $G$ be a finite group and $\7F$ an algebraically closed field whose characteristic does 
not divide $|G|$. There exists a strongly separable Frobenius algebra
$(\Gamma,m,\eta,\Delta,\ve)$ in $\2C=\Rep_\7F\,G$ such that 
\begin{enumerate}
\item $\DS\alpha\beta=|G|$. We normalize such that $\beta=1$.
\item $\Gamma$ is (isomorphic to) the left regular representation of $G$,
\item $\Gamma\otimes X\cong d(X)\Gamma\ \ \forall X$.
\item $\DS\dim\Hom_\2C(\11,\Gamma)=1$.
\end{enumerate}
If $\7F=\7C$, $\Rep\,G$ is a $*$-category and one can achieve $\Delta=m^*,\ve=\eta^*$.
\eprop

\brem 1. See also \cite{brug1} where a similar, but less symmetric, statement appears.

2. The proposition generalizes to finite dimensional Hopf algebras $H$, where the
categorical Frobenius algebra in $H\!-\!\Mod$ is strongly separable iff $H$ is semisimple
and cosemisimple, cf.\ \cite{mue09}.

3. Some of the structure survives for infinite compact groups and discrete quantum groups,
cf. \cite{mue13}.
\erem

\brem \label{rem-aut}
Given the monoid part of the above Frobenius algebra, one can obtain a fiber functor
$E: \2C\rarr\mathrm{Vect}_\7F$ as follows:
\bean E(X) &=& \Hom_\2C(\11,\Gamma\otimes X), \\
 E(s)\phi &=& s\otimes\id_X\mcirc\phi, \quad s:X\rarr Y, \ \phi\in E(X). 
\eean
The natural isomorphisms $d_{X,Y}: E(X)\otimes E(Y)\rarr E(X\otimes Y)$ are given by
\[ d_{X,Y}(\phi\boxtimes\psi)
   =m\otimes\id_X\otimes\id_Y\mcirc\id_\Gamma\otimes\phi\otimes\id_Y \mcirc\psi, \quad
   \phi\in E(X), \ \psi\in E(Y). \]
(Similarly, one can use the comonoid structure.) For the details, which are an immediate
generalization of \cite{del}, see \cite{mue13}. A similar construction is given in
\cite[Appendix C]{y1}. Defining  
\[ \Aut(\Gamma,m,\eta)\equiv \{g\in\End\,\Gamma\ | \ g\circ m=m\circ g\otimes g, \
   g\circ\eta=\eta\} \] 
it is easy to see that 
\[ g\mapsto (g_X), \quad g_X(\phi)=g\otimes\id_X\mcirc\phi, \ \ X\in\2C,\ \phi\in E(X) \]
defines a homomorphism $\Aut(\Gamma,m,\eta)\rarr\mathrm{Nat}^\otimes\,E=G$.
Appealing to the Yoneda lemma one verifies that this is a bijection, implying that 
$\Aut(\Gamma,m,\eta)$ is a group. This allows to recover $G$ from the
monoid structure on the regular representation without reference to the fiber functor
arising from the latter. This will turn out very useful in the sequel.
\erem

\brem
In fact, in \cite{mue13} a proof of Theorem \ref{theor2} will be given, whose first
step is to construct from a category $\2C$ (not necessarily finite) a monoid
$(\Gamma,m,\eta)$ (in $\mathrm{Ind}\,\2C$ if $\2C$ is infinite) such that 
$\Gamma\otimes X\cong d(X)\Gamma$ and $\dim\Hom(\11,\Gamma)=1$. One then obtains $G$
simply as the automorphism group of the monoid as above, the monoid of course turning out
to be the regular monoid of $G$. (This goes beyond the proof in \cite{del} that used a
monoid not satisfying the latter condition. This monoid is not the regular representation
and gives rise to a fiber functor into $\mathrm{Vect}_\7F$ only after a quotient
operation. Thus one cannot define $G$ as the automorphism group of the monoid.)
\erem

Even though the only monoids and Frobenius algebras considered in this paper are those
arising from regular representations as in Proposition \ref{prop2}, it is natural to give
the following considerations in larger generality.

\bdefprop \cite{par} Let $\2C$ be a strict tensor category and let $(\Gamma,m,\eta)$ be a
monoid in $\2C$. A $\Gamma$-module in $\2C$ is a pair $(X,\mu)$ where $X\in\2C$ and 
$\mu: \Gamma\otimes X\rarr X$ satisfies  
\[ \mu\mcirc\id_\Gamma\otimes\mu=\mu\mcirc m\otimes\id_X,
 \quad \quad \mu\mcirc\eta\otimes\id_X=\id_X. \]
The modules form a category $\Gamma\!-\!\Mod_\2C$ where 
$\Hom_{\Gamma\!-\!\Mod}((X,\mu),(Y,\lambda))=\{s: X\rarr Y\ |\ s\circ\mu=\lambda\circ\id_\Gamma\otimes s\}$. 
If $\2C$ is braided and has coequalizers, $\otimes$ preserves coequalizers, and 
$(\Gamma,m,\eta)$ is commutative then $\Gamma\!-\!\Mod$ is a tensor category with 
$(X,\mu)\otimes (Y,\eta)=\mathrm{coeq}(\alpha,\beta)$, where 
$\alpha,\beta: \Gamma\otimes X\otimes Y\rarr X\otimes Y$ are given by
\[ \alpha=\mu\otimes\id_Y, \quad\quad
   \beta=\id_X\otimes\eta\mcirc c_{\Gamma,X}\otimes\id_Y. \]
The full subcategory $\Gamma\!-\!\Mod_\2C^0\subset\Gamma\!-\!\Mod_\2C$ consisting of the objects
$(X,\mu)$ satisfying $\mu\circ c_{X,\Gamma}\circ c_{\Gamma,X}=\mu$ is monoidal and
braided.
\edefprop

\brem
1. The above definition and facts are due to Pareigis \cite{par} and were rediscovered in
\cite{ko}. The special case where $\Gamma\in Z_2(\2C)$, implying
$\Gamma\!-\!\Mod^0_\2C=\Gamma\!-\!\Mod_\2C$, was considered in \cite{brug1}.

2. Note that the coequalizers are unique only up to isomorphism, thus some care
is required in the definition of the associativity constraint of $\Gamma\!-\!\Mod_\2C$. In
\cite{par} this is handled by showing that $\Gamma\!-\!\Mod_\2C$ is equivalent (as a
category) to a full subcategory of the category of $M-M$ bimodules in $\2C$. For the
latter the associativity constraint had been constructed in B.\ Pageigis: Non-additive
ring and module theory V. Algebra Berichte {\bf 40}, 1980.

3. We will exclusively consider semisimple categories with duals. In such categories,
coequalizers exist and are preserved by $\otimes$.
\erem

Recall that the dimension of a finite TC is the sum over the squared dimensions of its
simple objects, cf.\ e.g.\ \cite{bw, mue09}.

\bprop \label{prop-dim}
Let $\2C$ be a finite BTC and let $(\Gamma,m,\eta,\Delta,\ve)$ be a strongly
separable Frobenius algebra in $\2C$ satisfying $\dim\Hom(\11,\Gamma)=1$. Then
$\Gamma\!-\!\Mod_\2C$ is a semisimple $\7F$-linear spherical tensor category with
$\End_\Gamma\11=\7F\id_\11$, and 
\[ \dim \Gamma\!-\!\Mod_\2C = (\dim\Gamma)^{-1}\,\dim\2C. \]
\eprop

\prf The free module functor 
$F: \2C\rarr\Gamma\!-\!\Mod,\ X\mapsto(\Gamma\otimes X,m\otimes\id_X)$ is a left adjoint
of the forgetful functor $G: \Gamma\!-\!\Mod\rarr\2C,\ (X,\mu)\mapsto X$, cf.\
\cite{brug1,ko}. $F$ is monoidal, implying $F(\11)\cong\11$ and $d(F(X))=d(X)$.
The tensor unit of $\Gamma\!-\!\Mod_\2C$ being $(\Gamma,m)$ we have 
$\End_\Gamma\11=\Hom_\Gamma(F(\11),(\Gamma,m))\cong\Hom(\11,\Gamma)$, implying
$\End_\Gamma\11=\7F\id_\11$. As a rigid ribbon category, $\2C$ is spherical and so is 
$\Gamma\!-\!\Mod_\2C$ \cite{mue09}, allowing us to talk of dimensions of
objects irrespective of whether $\Gamma\!-\!\Mod_\2C$ is braided.
Semisimplicity is proven as in \cite{brug1,ko}; it is here that the Frobenius
structure is used, cf.\ also \cite{mue09}. The fact $GF(X)=\Gamma\otimes X$  
together with $d(F(X))=d(X)$ and additivity of $F$ and $G$ implies
$d(G(Y))=d(\Gamma)d(Y)$. Let now $\{X_i\in\2C\}$ and $\{Y_j\in\Gamma\!-\!\Mod_\2C\}$ be
complete sets of simple objects in $\2C$ and $\Gamma\!-\!\Mod_\2C$, respectively. The
computation 
\bean \dim\2C &=& \sum_i d(X_i)^2 \ =\ \sum d(X_i)d(F(X_i)) \\
  &=& \sum_i\sum_j d(X_i)d(Y_j)\dim\Hom(F(X_i),Y_j) \\
  &=& \sum_i\sum_j d(X_i)d(Y_j)\dim\Hom(X_i,G(Y_j)) \\
  &=& \sum_j d(Y_j)d(G(Y_j))\ =\ d(\Gamma)\sum_j d(Y_j)^2 \\
  &=& d(\Gamma)\dim\Gamma\!-\!\Mod_\2C \eean
completes the proof.
\qed

\brem A similar result is proven in \cite{brug1} where $\Gamma\in Z_2(\2C)$, implying
$\Gamma\!-\!\Mod$ to be braided, is assumed. The present very simple proof shows that such
an assumption is not needed.
\erem

While the category $\Gamma\!-\!\Mod_\2C$ considered above is conceptually very natural, there
is an alternative description which occasionally is more convenient. The point is that the
tensor product of $\Gamma\!-\!\Mod_\2C$, while entirely analogous to that in $R-\Mod$, is not
very convenient to work with.

\bdefprop \label{defin_galois}
Let $\2C$ be a strict BTC and $(\Gamma,m,\eta,\Delta,\ve)$ a strongly separable Frobenius
algebra in $\2C$. Then the following defines a tensor category $\tilde{\2C}_\Gamma$.
\begin{itemize}
\item $\Obj\tilde{\2C}_\Gamma=\Obj\2C$.
\item $X\tilde{\otimes}Y=X\otimes Y$.
\item $\Hom_{\tilde{\2C}_\Gamma}(X,Y)=\Hom_\2C(\Gamma\otimes X,Y)$.
\item Let $s\in\Hom_{\tilde{\2C}_\Gamma}(X,Y)=\Hom_\2C(\Gamma\otimes X,Y)$ and
$t\in\Hom_{\tilde{\2C}_\Gamma}(Y,Z)=\Hom_\2C(\Gamma\otimes Y,Z)$. Then
$t\tilde{\circ} s=t\circ\id_\Gamma\otimes s\circ\Delta\otimes X$ in
$\Hom_{\tilde{\2C}_\Gamma}(X,Z)=\Hom_\2C(\Gamma\otimes X,Z)$.
\item Let $s\in\Hom_{\tilde{\2C}_\Gamma}(X,Y)=\Hom_\2C(\Gamma\otimes X,Y)$ and
$t\in\Hom_{\tilde{\2C}_\Gamma}(Z,T)=\Hom_\2C(\Gamma\otimes Z,T)$. Then
$s\tilde{\otimes}t=s\otimes t\mcirc\id_\Gamma\otimes c_{\Gamma,X}\otimes\id_Z\mcirc\Delta\otimes\id_X\otimes\id_Z$  
in $\Hom_{\tilde{\2C}_\Gamma}(X\otimes Z,Y\otimes T)=\Hom_\2C(\Gamma\otimes X\otimes Z,Y\otimes T)$.
\end{itemize}
The canonical completion $\hat{\2C}_\Gamma=\tilde{\2C}_\Gamma^p$ of $\tilde{\2C}_\Gamma$
to a category with splitting idempotents is semisimple. (Recall that
$\Obj\,\hat{\2C}_\Gamma=\{(X,p),\ X\in\Obj\tilde{\2C}_\Gamma,\, p=p^2\in\End_{\tilde{\2C}_\Gamma}X\}$
etc. Instead $(X,\id_X)\in\hat{\2C}_\Gamma$ we simply write $X$.)
If $\2C$ is a $*$-category and $\Delta=m^*, \ve=\eta^*$ then
$\tilde{\2C}_\Gamma, \hat{\2C}_\Gamma$ are $*$-categories.
The functor $\iota: \2C\rarr\tilde{\2C}_\Gamma$ given by 
$X\mapsto X, \ s\mapsto \ve\otimes s$ is monoidal and faithful. The composite of $\iota$
with the full embedding $\tilde{\2C}_\Gamma\rarr\hat{\2C}_\Gamma$ is also denoted by
$\iota$.
\edefprop

\begin{figure}
\[ \ba{ccc} \begin{tangle}
\step[1.5]\object{Z} \\
\hh\step\hstep\id \\
\step\frabox{Z} \\
\hh\step\id\step\id\step[.3]\obj{Y} \\
\hh\step\id\step\frabox{t} \\
\hh\step[.2]\obj{\Gamma}\step[.8]\id\step[.3]\obj{\Gamma}\step[.7]\id\step\id \\
\hh\hstep\obj{\Delta}\hstep\cu\step\id\\
\hh\step[1.5]\id\step[1.5]\id \\
\step[1.5]\object{\Gamma}\step[1.5]\object{X}
\end{tangle} 
& \quad\quad &
\begin{tangle}
\step[1.5]\object{Y}\step[3]\object{T} \\
\hh\step[1.5]\id\step[3]\id \\
\step\frabox{s}\Step\frabox{t} \\
\step[.2]\obj{\Gamma}\step[.8]\id\hstep\obj{\Gamma}\hstep\xx\step\id \\
\hh\hstep\obj{\Delta}\hstep\cu\Step\id\step\id \\
\hh\step[1.5]\id\step[2.5]\id\step\id \\
\step[1.5]\object{\Gamma}\step[2.5]\object{X}\step\object{Z}
\end{tangle} \\
\hbox{Composition} && \hbox{Tensor Product} \ea \]
\caption{Composition and Tensor product of arrows in $\2C\rtimes_0\2S$}
\label{fig1}\end{figure}

\prf That $\tilde{\2C}_\Gamma$ and therefore $\hat{\2C}_\Gamma$ is a $\7F$-linear strict
tensor category is almost obvious: One only needs to show associativity of 
$\tilde{\circ}, \tilde{\otimes}$ on the morphisms and the interchange law, which is left
to the reader. The discussion of the $*$-operation on $\tilde{\2C}_\Gamma,\hat{\2C}_\Gamma$
and of semisimplicity of $\hat{\2C}_\Gamma$ is the same as in \cite{mue06,mue09}, to which
we refer for details.
\qed

\bprop \label{prop_equiv}
Let $\2C$ and $(\Gamma,m,\eta,\Delta,\ve)$ be as before. Then there exists a monoidal
equivalence $K: \hat{\2C}_\Gamma\rarr\Gamma\!-\!\Mod_\2C$ such that $K\circ\iota\cong F$ as
tensor functors.
\eprop

\prf We define $K_0: \tilde{\2C}_\Gamma\rarr\Gamma\!-\!\Mod_\2C$ by $K_0(X)=F(X)$, and for
$s\in\Hom_{\tilde{\2C}_\Gamma}(X,Y)=\Hom(\Gamma\otimes X,Y)$ we put 
$K_0(s)=\id_\Gamma\otimes s\mcirc\Delta\otimes\id_X\in K_0(s)\in\Hom_\Gamma(F(X),F(Y))$.
The map $s\mapsto K_0(s)$ has inverse $t\mapsto\ve\otimes\id_Y\mcirc t$. Direct
computations show $K_0(s)\bullet K_0(t)=K_0(s\bullet t)$ for
$\bullet\in\{\circ,\otimes\}$, thus $K$ is a full and faithful tensor functor and
satisfies $K_0\circ\iota=F$. Since $\Gamma\!-\!\Mod$ has splitting idempotents, 
$K_0: \tilde{\2C}_\Gamma\rarr\Gamma\!-\!\Mod$ extends to 
$K: \hat{\2C}_\Gamma\rarr\Gamma\!-\!\Mod$, uniquely up to natural isomorphism. Since every
object of $\Gamma\!-\!\Mod$ is a retract of an object of the form $K_0(X)=F(X)$, $K$ is 
essentially surjective, thus an equivalence. 
\qed


\section{Braided Crossed G-Category from Galois Extensions} \label{sec3}
\subsection{Definition of $\2C\rtimes\2S$ and Basic Properties}
In the rest of the paper we assume $\7F$ to satisfy the assumptions of Theorem
\ref{theor2} and work exclusively with the category $\hat{\2C}_\Gamma$. Furthermore, $\2C$
will be a BTC, not necessarily finite, and $\2S\subset\2C$ will be a finite admissible
full sub-STC. 

\blemma \label{lemma_fixp}
Let $G$ be such that $\2S\simeq\Rep\,G$ and let $(\Gamma,\ldots)$ be the corresponding
commutative Frobenius algebra in $\2C$. We write $p_0=\eta\circ\ve\in\End\,\Gamma$ and
recall that $G=\Aut(\Gamma,m,\eta)$. For $s\in\Hom(\Gamma\otimes X,Y)$ the following are
equivalent:
\begin{itemize}
\item[(i)] $s\mcirc g\otimes\id_X=s$ for all $g\in G$.
\item[(ii)] $s\mcirc p_0\otimes\id_X=s$.
\end{itemize}
\elemma

\prf (ii)$\impl$(i): Obvious consequence of $\ve\circ g=\ve\ \forall g\in G$.

(i)$\impl$(ii): If $\hat{G}$ denotes the set of iso-classes of irreps $\pi_i$ of $G$ and 
$d_i$ is the dimension of $\pi_i$, we have 
$\End\,\Gamma\cong\oplus_{i\in\hat{G}} M_{d_i}(\7F)$ and 
$G\ni g=\oplus_{i\in\hat{G}}\,\pi_i(g)$. Whenever $\hat{G}\ni i\ne 0$ there exists 
$g\in G$ such that $\pi_i(g)\ne\id$. If $p_i$ is the unit of $M_{d_i}(\7F)$, (i) therefore
implies $s\mcirc p_i\otimes\id_X=0\ \forall i\ne 0$, and we conclude
$s=\sum_i s\mcirc p_i\otimes\id_X=s\mcirc p_0\otimes\id_X$.
\qed

\bdefin Let $\2C$ be a strict BTC and $\2S\subset\2C$ a finite full sub-STC. Let
$(\Gamma,\ldots)$ be the Frobenius algebra in $\2C$ arising from Theorem \ref{theor2} and
Proposition \ref{prop2}. Then we write $\2C\rtimes_0\2S:=\tilde{\2C}_\Gamma$ and
$\2C\rtimes\2S:=\hat{\2C}_\Gamma$. 
\edefin

For the sake of legibility we will continue to write $\tilde{\2C}, \hat{\2C}$ rather than
$\2C\rtimes_0\2S, \2C\rtimes\2S$ in many places, in particular subscripts.

\bprop \label{prop-cs}
$\2C\rtimes_0\2S$ and $\2C\rtimes\2S$ are strict spherical tensor categories and
$\2C\rtimes\2S$ is semisimple. If $\2C$ is a $*$-category then $\2C\rtimes_0\2S$ and
$\2C\rtimes\2S$ have a $*$-structure extending that of $\2C$. There exists a canonical
tensor functor $\iota: \2C\rarr\2C\rtimes\2S$ which is faithful and injective on the
objects, thus an inclusion. The group $G=\Aut(\Gamma,m,\eta)$ acts on $\2C\rtimes\2S$ via
$\gamma_g(s)=s\mcirc g^{-1}\otimes\id_X$ for  
$s\in\Hom_{\2C\rtimes\2S}(X,Y)=\Hom(\Gamma\otimes X,Y)$ and
$\gamma_g((X,p))=(X,\gamma_g(p))$. We have $(\2C\rtimes_0\2S)^G\cong\2C$ 
and $(\2C\rtimes\2S)^G\simeq\2C$. If $\2C$ is finite then 
$\dim\2C\rtimes\2S=\dim\2C/|G|=\dim\2C/\dim\2S$.
\eprop

\prf The first set of statements is obvious. Clearly, $g\mapsto\gamma_g$ is a homomorphism
and $\gamma_g$ is invertible. Now $\gamma_g(s\bullet t)=\gamma_g(s)\bullet\gamma_g(t)$ for
$\bullet\in\{\circ,\otimes\}$ follows from $\Delta\circ g=g\otimes g\circ\Delta$. Lemma 
\ref{lemma_fixp} amounts to $(\2C\rtimes_0\2S)^G=\iota(\2C)\cong\2C$, and
$(\2C\rtimes\2S)^G=\ol{\iota(\2C)}^p\cong\2C^p\simeq\2C$. The dimension formula follows
from Propositions \ref{prop-dim} and \ref{prop_equiv}.
\qed

\brem 1. Here and in the sequel, $\2D^G\subset\2D$ denotes the subcategory consisting of the
objects and morphisms that are strictly fixed by the action of $G$. In our strict context
this is the right notion, but it presumably needs to be generalized if one works with a less
strict notion of G-categories. 

2. For definition of $\2C\rtimes\2S$ given above for finite $\2S$ is equivalent to the one
in \cite{mue06}. Thus Proposition \ref{prop_equiv} proves the equivalence of the approaches 
to Galois extensions and modularization of braided tensor categories given by the author
\cite{mue06} and A. Brugui\`{e}res \cite{brug1}. While both definitions are equally
involved, $\Gamma\!-\!\Mod_\2C$ may be more natural, yet $\hat{\2C}_\Gamma$ has some
advantages. On the one hand, the tensor product of $\hat{\2C}_\Gamma$ is canonical, i.e.\
involving no choices, and strict, making it slightly more convenient to work with. On the
other hand, the relationship between the categorical constructions and (algebraic) quantum
field theory, cf.\ the next section, is very easy to establish for $\2C\rtimes\2S$.

3. When $\2S$ is infinite the definition of $\2C\rtimes\2S$ must be changed. While there
still is a monoid structure on the regular representation $\Gamma$ \cite{mue13}, the
latter lives in a larger category $\mathrm{Ind}\,\2S$ and is no more a Frobenius
algebra. Thus the proof of semisimplicity also changes. The somewhat pedestrian
definition of $\2C\rtimes\2S$ in \cite{mue06} works also for infinite $\2S$.

4. Constructions similar to the one above have been given in \cite{brug1, y1,y2,fs}.
\erem

The following is due to Brugui\`eres \cite{brug1}, who proved it for the category of 
$(\Gamma,m,\eta)$-modules.

\btheor \label{th-univ}
Let $\2S\subset\2C$ be as before. The tensor functor $\iota: \2C\rarr\2C\rtimes\2S$ has
the following universal property:
\begin{enumerate}
\item $\iota$ is faithful and for every simple object $Y\in\2C\rtimes\2S$ there exists
$X\in\2C$ such that $Y$ is a direct summand of $Y\prec \iota(X)$. 
\item For every $X\in\2S$ we have $\iota(X)\cong d(X)\11$ in $\2C\rtimes\2S$.
\item If $\2D$ is semisimple and $\iota': \2C\rarr\2D$ satisfies 1-2 then there exists a
faithful tensor functor $\iota'': \2C\rtimes\2S\rarr\2D$, unique up to monoidal natural
isomorphism, such that $\iota'=\iota''\circ\iota$.
\end{enumerate}
\etheor

\prf 1. Obvious by construction. 2. It is sufficient to show this for $X\in\2S$ simple.
We have $\Hom_{\hat{\2C}}(\11,\iota(X))=\Hom_\2C(\Gamma,X)$, and $\Gamma\cong\oplus_i d(X_i)X_i$ 
implies $\dim\Hom_{\hat{\2C}}(\11,\iota(X))=d(X)$. Thus $\iota(X)\cong d(X)\11\oplus X'$ and 
$\End\,\iota(X)\cong M_{d(X)}\oplus N$. Now, 
\bean \dim\End_{\hat{\2C}}\iota(X) &=& \dim\Hom_\2C(\Gamma\otimes X,X)
  =\dim\Hom_\2C(d(X)\Gamma,X)  \\ &=& d(X)\dim\Hom_\2C(\Gamma,X)=d(X)^2, \eean
thus $N=0$ and $\iota(X)\cong d(X)\11$.

3. This follows from the corresponding statement in \cite{brug1} and Proposition \ref{prop_equiv}.
(We omit the direct proof for reasons of space.)
\qed

The considerations in the remainder of this section concern the decomposition of
$\iota(X)\in\2C\rtimes\2S$ for simple $X\in\2C$, complementing the results in
\cite[Section 4.1]{mue06}, and will not be used in the rest of the paper. 

\bdefin For $X,Y\in\2C$ we write $X\sim Y$ iff $\Hom_\2C(\Gamma\otimes X,Y)\ne\{0\}$. 
\edefin

\btheor 
Restricted to simple objects, the relation $\sim$ is an equivalence relation.
Let $X,Y\in\2C$ be simple. If $X\not\sim Y$ then $\iota(X), \iota(Y)$ are disjoint, to wit 
$\iota(X),\iota(Y)$ have no isomorphic subobjects. For every equivalence class $\sigma$ there
exist a finite set $\2I_\sigma$, mutually non-isomorphic simple objects
$Z_i\in\2C\rtimes\2S,\ i\in\2I_\sigma$ and natural numbers $N_X, X\in\sigma$ such that
\[ \iota(X)\cong N_X\ \bigoplus_{i\in\2I_\sigma} Z_i \quad\forall X\in\sigma. \]
\etheor

\prf For all $X,Y$ we have $X\sim X$ (since $\11\prec\Gamma$) and
$X\sim Y\Leftrightarrow Y\sim X$ (since $\ol{\Gamma}\cong\Gamma$). Let $X,Y,Z $ be simple
and $X\sim Y\sim Z$. $\Hom(\Gamma\otimes X,Y)\ne\{0\}$ implies $Y\prec\Gamma\otimes X$,
i.e.\ $Y$ is a direct summand of $\Gamma\otimes X$. Similarly, $X\prec\Gamma\otimes Y$,
$Y\prec\Gamma\otimes Z$, $Z\prec\Gamma\otimes Y$. Thus
$X\prec\Gamma\otimes Y\prec\Gamma\otimes\Gamma\otimes Z\cong|G|\,\Gamma\otimes Z$, where we used
$\Gamma^2\cong|G|\Gamma$. Therefore $\Hom(X,\Gamma\otimes Z)\ne\{0\}$, thus $X\sim Z$, and
$\sim$ is transitive. In view of 
$\Hom(\Gamma\otimes X,Y)=\Hom_{\hat{\2C}}(\iota(X),\iota(Y))$ it is clear that $X\not\sim
Y$ implies disjointness.

Let $X,Y\in\2C$ be simple such that $X\sim Y$ and let $Z_1\prec\iota(X)$ be simple. 
Together with $\iota(X)\prec\iota(\Gamma)\iota(Y)$ this implies 
$Z_1\prec\iota(\Gamma)\iota(Y)\prec|G|\iota(Y)$, where we used $\iota(\Gamma)\cong|G|\11$.  
Since $Z_1$ is simple, we have $Z_1\prec\iota(Y)$. Thus every simple $Z_1\in\2C\rtimes\2S$
contained in $\iota(X)$ is also contained in $\iota(Y)$, provided $X\sim Y$. We conclude
that $X\sim Y$ implies that $\iota(X)$ and $\iota(Y)$ contain the same simple
summands. The rest follows from the fact \cite[Proposition 4.2]{mue06} that, for every
simple $X\in\2C$, the simple summands of $\iota(X)\in\2C\rtimes\2S$ appear with the same
multiplicity $N_X$.
\qed

\brem If $G$ is abelian, corresponding to all simple objects in $\2S$ being invertible, 
we have $X\sim Y$ iff there exists $Z\in\2S$ such that $X\cong Z\otimes Y$. As a consequence,
$X\sim Y$ implies $\iota(X)\cong\iota(Y)$ and $N_X=N_Y$. 
Since in the general case $X\sim Y$ does not imply that $X,Y$ have the same dimensions the
above result, according to which $\iota(X),\iota(Y)$ have the same simple summands, clearly
is the best one can hope for.

In the abelian case, the structure of $\End_{\hat{\2C}}\iota(X)$ can be clarified quite
explicitly, cf.\ \cite[Sect.\ 5.1]{mue06}. Presently there is no analogous result in the
general, non-abelian case. 
\erem


\subsection{$\2C\rtimes\2S$ as Braided Crossed G-category}
Let $c$ be the braiding of $\2C$. For $X,Y\in\2C$ it is clear that $\iota(c_{X,Y})$ is an
isomorphism $\iota(X)\iota(Y)\rarr\iota(Y)\iota(X)$ satisfying the braid
equations. Whether this gives rise to a braiding of $\2C\rtimes_0\2S$ (and therefore of
$\2C\rtimes\2S$) depends on whether or not $\iota(c)$ is natural w.r.t.\ the larger
hom-sets of $\2C\rtimes_0\2S$. For one variable we in fact have:

\blemma \label{lem-nat1}
Let $X,Y,Z\in\2C$ and 
$s\in\Hom_{\hat{\2C}}(X,Y)=\Hom_\2C(\Gamma\otimes X,Y)$. Then 
\[ \iota(c_{Y,Z})\ \hat{\circ}\ s\hat{\otimes}\id_Z=\id_Z\hat{\otimes} s\ \hat{\mcirc}\ 
   \iota(c_{X,Z}) \]
holds in $\2C\rtimes_0\2S$.
\elemma

\prf In view of Definition \ref{defin_galois}, the two sides of the desired equation are
represented by the following morphisms in $\2C$:
\[ \begin{array}{c|c}
   \Hom_{\hat{\2C}}(X\otimes Z,Z\otimes Y) & \Hom(\Gamma\otimes X\otimes Z,Z\otimes Y) \\ \hline
 \iota(c_{Y,Z})\ \hat{\circ}\ s\hat{\otimes}\id_Z & c_{Y,Z}\mcirc s\otimes\id_Z \\
  \id_Z\hat{\otimes} s\ \hat{\circ}\ \iota(c_{X,Z}) & 
       \id_Z\otimes s\mcirc c_{\Gamma,Z}\otimes\id_X\mcirc\id_\Gamma\otimes c_{X,Z}
\end{array}\]
A trivial computation in $\2C$ shows that the expressions on the right hand side coincide. 
\qed

As shown in \cite{mue06}, naturality of $c$ w.r.t.\ the second variable holds iff
$\2S\subset Z_2(\2C)$, which is the case iff $\Gamma\in Z_2(\2C)$. Here
$Z_2(\2C)\subset\2C$ is the full subcategory of objects $X$ satisfying 
$c_{X,Y}\circ c_{Y,X}=\id_{YX}$ for all $Y\in\2C$, called central in \cite{mue11} and
transparent in \cite{brug1}. In order to understand the general case $\2S\not\subset
Z_2(\2C)$ we need some preliminary considerations. 

\blemma \label{lem-nat2}
Let $X,Y\in\2C, \ Z\in\2C\cap\2S'$ and $s\in\Hom_{\hat{\2C}}(X,Y)=\Hom_\2C(\Gamma\otimes X,Y)$. Then 
\[ \iota(c_{Z,Y})\ \hat{\circ}\ \id_Z\hat{\otimes} s= s\hat{\otimes}\id_Z
   \ \hat{\circ}\ \iota(c_{Z,X}). \]
\elemma

\prf As above we have
\[ \begin{array}{c|c}
   \Hom_{\hat{\2C}}(Z\otimes X,Y\otimes Z) & \Hom(\Gamma\otimes Z\otimes X,Y\otimes Z) \\ \hline
   \iota(c_{Z,Y})\ \hat{\circ}\ \id_Z\hat{\otimes} s & 
      c_{Z,Y}\mcirc \id_Z\otimes s\mcirc c_{\Gamma,Z}\otimes\id_X \\
   s\hat{\otimes}\id_Z \ \hat{\circ}\ \iota(c_{Z,X}) & 
        s\otimes\id_Z\mcirc\id_\Gamma\otimes c_{Z,X}
\end{array}\]
Now we find
\bean \lefteqn{ c_{Z,Y}\mcirc \id_Z\otimes s\mcirc c_{\Gamma,Z}\otimes\id_X } \\
  &&=   s\otimes\id_Z\mcirc \id_\Gamma\otimes c_{Z,X}\mcirc
   (c_{Z,\Gamma}\mcirc c_{\Gamma,Z})\otimes\id_X. \eean
For arbitrary $Z\in\2C$ this will not coincide with
$s\otimes\id_Z\mcirc\id_\Gamma\otimes c_{Z,X}$, but for $Z\in\2C\cap\2S'$ it does since 
$\Gamma\in\2S$, implying $c_{Z,\Gamma}\mcirc c_{\Gamma,Z}=\id$.
\qed
\\

Let $X\in\2C$ and $p\in\End_{\hat{\2C}}(X)$ a minimal idempotent, thus
$X_1=(X,p)\in\2C\rtimes\2S$ is simple. Let $v: X_1\rarr(X,p),\ v': (X,p)\rarr X_1$
satisfy $v'\circ v=\id_{X_1}$, $v\circ v'=p$ and consider
\begin{equation} \label{eq-del0}
\id_{\Gamma}\hat{\otimes}v\ \hat{\mcirc}\ \iota(c_{X,\Gamma}\circ c_{\Gamma,X})\
   \hat{\circ}\ \id_{\Gamma}\hat{\otimes}v' \ \in\End_{\hat{\2C}}(\Gamma\otimes X_1).
\end{equation}
In view of $\Gamma\in\2S\subset\2C\cap\2S'$, the preceding lemmas imply that
(\ref{eq-del0}) equals 
\[ \id_{\Gamma}\hat{\otimes}p\ \hat{\mcirc}\ \iota(c_{X,\Gamma}\circ c_{\Gamma,X})
   = \iota(c_{X,\Gamma}\circ c_{\Gamma,X})\ \hat{\circ}\ \id_{\Gamma}\hat{\otimes}p, \]
which in particular implies that (\ref{eq-del0}) is invertible, thus is in
$\Aut_{\hat{\2C}}(\Gamma\otimes X_1)$. The inverse is given by 
\[ \del''X_1\ := \ \id_\Gamma\hat{\otimes}v\ \hat{\mcirc}\ 
    \iota(\tilde{c}(X,\Gamma)\circ\tilde{c}(\Gamma,X))\ \hat{\circ}\
   \id_\Gamma\hat{\otimes}v' \ \in\Aut_{\hat{\2C}}(\Gamma\otimes X_1), \]
where $\tilde{c}(X,Y)=c_{Y,X}^{-1}$.
Since $X_1$ is simple and $\iota(\Gamma)\cong |G|\,\id_\11$ we have
\begin{equation} \label{e1}
  \del''X_1\ = \ \del' X_1\,\hat{\otimes}\,\id_{X_1}, 
\end{equation}
where $\del' X_1\in\Aut\iota(\Gamma)\cong M_{|G|}(\7F)$. This equation, which lives in
$\hat{\2C}$, corresponds to
\[ \iota(\tilde{c}(X,\Gamma)\circ\tilde{c}(\Gamma,X))\ \hat{\circ}\
   \id_\Gamma\hat{\otimes}p=\del'X_1\hat{\otimes}p \]
in $\tilde{\2C}$ and to 
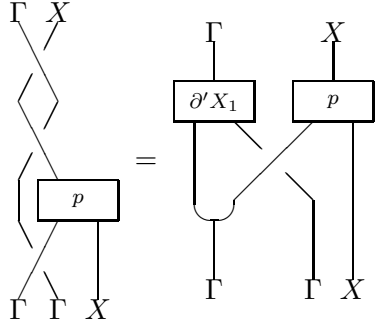
\[
\begin{tangle}
\object{\Gamma}\step\object{X}\\
\hx\\
\hx\\
\hh\id\step\frabox{p}\\
\hxx\step\id\\
\object{\Gamma}\step\object{\Gamma}\step\object{X}
\end{tangle}
\quad=\quad
\begin{tangle}
\hstep\object{\Gamma}\step[3]\object{X}\\
\hh\hstep\id\step\Step\id\\
\frabox{\del'X_1}\Step\frabox{p}\\
\id\step\xx\step\id\\
\hcu\Step\id\step\id\\
\hstep\object{\Gamma}\step[2.5]\object{\Gamma}\step\object{X}
\end{tangle}
\]
in $\2C$. Composing with $\Delta\otimes\id_X$ and using cocommutativity
$c_{\Gamma,\Gamma}\circ\Delta=\Delta$ we obtain
\begin{equation} 
\begin{tangle}
\object{\Gamma}\step\object{X}\\
\hx\\
\hx\\
\hh\id\step\frabox{p}\\
\hh\step[-.5]\obj{\Delta}\step[.5]\hcu\step\id\\
\hh\hstep\id\step[1.5]\id\\
\hstep\object{\Gamma}\step[1.5]\object{X}
\end{tangle}
\quad\quad=\quad\quad
\begin{tangle}
\hstep\object{\Gamma}\step[3]\object{X}\\
\hh\hstep\id\step\Step\id\\
\frabox{\del'X_1}\step[1.5]\frabox{p}\\
\hh\hcu\step[1.5]\id\step\id\\
\hstep\cu\step\id\\
\step[1.5]\object{\Gamma}\Step\object{X}
\end{tangle}
\quad\quad=\quad\quad
\begin{tangle}
\hstep\object{\Gamma}\step[2.5]\object{X}\\
\hh\hstep\id\step\step[1.5]\id\\
\frabox{\del X_1}\step[1.5]\frabox{p}\\
\hstep\cu\step\id\\
\step[1.5]\object{\Gamma}\Step\object{X}
\end{tangle}
\label{eq-twist}
\end{equation}
where we have defined
\[ \del X_1:=\del'X_1\mcirc\Delta\ \in\End_\2C(\Gamma). \]
Before we elucidate the significance of (\ref{eq-twist}) we derive an explicit formula for 
$\del(X,p)$. In view of (\ref{e1}) it is clear that
\bean \del'X_1 &=& d(X_1)^{-1} \,(\id_\Gamma\otimes Tr_{X_1})\del''X_1 \nn\\
  &=& d(X_1)^{-1}\,(\id_\Gamma\otimes Tr_X) \,
  [\iota(\tilde{c}(X,\Gamma)\circ\tilde{c}(\Gamma,X))\,\hat{\circ}\,\id_\Gamma\hat{\otimes}p].
\eean
We have $\del'X_1\in\End_{\hat{\2C}}\iota(\Gamma)$, and computation shows that
$\del'X_1\in\End_{\hat{\2C}}\iota(\Gamma)$ is represented by 
\[ d(X_1)^{-1} \, \id_\Gamma\otimes\ol{\ve}_X
   \mcirc(\tilde{c}(X,\Gamma)\circ\tilde{c}(\Gamma,X))\otimes\id_{\ol{X}}
   \mcirc\id_\Gamma\otimes p\otimes_{\ol{X}}\mcirc c_{\Gamma,\Gamma}\otimes\ve_X 
\]
in $\Hom_\2C(\Gamma\Gamma,\Gamma)$. Furthermore,
\bea d(X_1) &=& Tr_{X_1}(\id_{X_1})=Tr_{\iota(X)}(p)=Tr_X(p\mcirc \eta\otimes\id_X) \nn\\
  &=& \ol{\ve}_X\mcirc p\otimes\id_{\ol{X}}\mcirc \eta\otimes \ve_X. \label{eq-dim}\eea
For $\del X_1=\del'X_1\mcirc\Delta$ we thus obtain
\[ \del X_1= d(X_1)^{-1}\,\id_\Gamma\otimes\ol{\ve}_X
   \mcirc(\tilde{c}(X,\Gamma)\circ\tilde{c}(\Gamma,X)) 
   \otimes\id_{\ol{X}} \mcirc\id_\Gamma\otimes p\otimes\id_{\ol{X}}\mcirc 
   \Delta\otimes\ve_X, \]
where we have used the cocommutativity $c_{\Gamma,\Gamma}\circ\Delta=\Delta$.
In diagrammatic form:
\begin{equation} \label{e3}
\del X_1=\quad
\left( \ \ \
\begin{tangle}
\hh\step\coev\\
\hh\frabox{p}\step\id\\
\obj{\eta}\counit\step\hev\\
\end{tangle}
\ \ \right)^{-1}
\cdot\quad
\begin{tangle}
\object{\Gamma}\\
\hh\id\Step\coev\\
\x\step\id\obj{\ol{X}}\\
\x\step\id\\
\hh\id\step\frabox{p}\step\id\\
\hh\step[-.5]\obj{\Delta}\hstep\hcu\hstep\obj{X}\hstep\hev\\
\hh\hstep\id\\
\hstep\object{\Gamma}
\end{tangle}
\end{equation}

By definition, $\del(X,p)\in\End_\2C(\Gamma)$. In fact we have a much stronger result.

\bprop \label{prop_grade}
Let $(X,p)\in\2C\rtimes\2S$ be simple. Then $\del(X,p)\in G=\Aut(\Gamma,m,\eta)$.
\eprop

\prf Since $\dim\Hom(\11,\Gamma)=1$ we have $\del(X,p)\circ\eta=c\eta$ and
$\ve\circ\del(X,p)=c\ve$, where $c=\ve\circ\del(X,p)\circ\eta$. Thus
\[ c=d(X_1)^{-1}\ol{\ve}_X\mcirc p\otimes\id_{\ol{X}}\mcirc\eta\otimes\ve_X
\]
and comparison with (\ref{eq-dim}) shows $c=1$, thus $\del(X,p)\circ\eta=\eta$.
Next, we compute
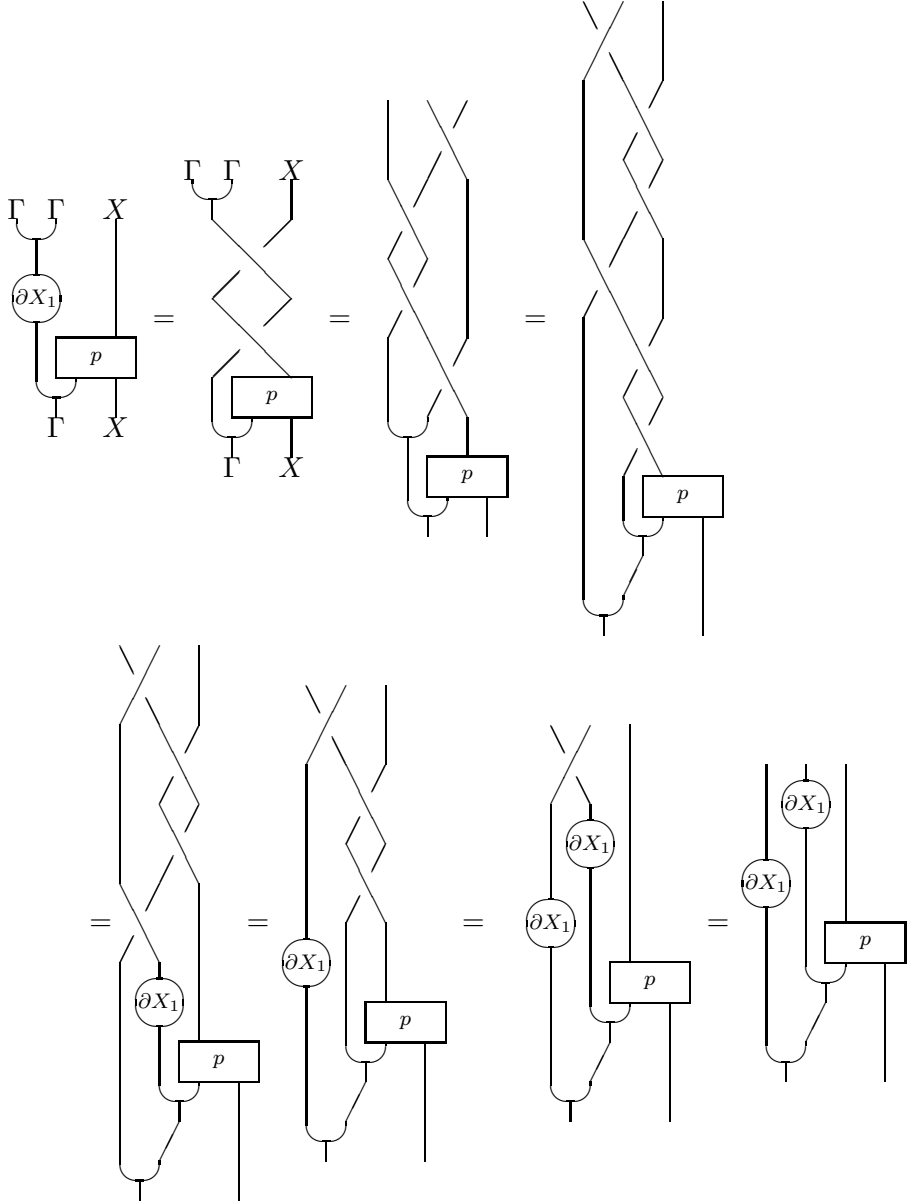
\bean \lefteqn{ 
\begin{tangle}
\step[-.5]\object{\Gamma}\step\object{\Gamma}\step[1.5]\object{X}\\
\hh\step[-.5]\hcu\step[1.5]\id\\
\O{\del X_1}\Step\id\\
\hh\id\step\frabox{p}\\
\hh\hcu\step\id\\
\hstep\object{\Gamma}\step[1.5]\object{X}
\end{tangle}
\quad=\quad
\begin{tangle}
\step[-.5]\object{\Gamma}\step\object{\Gamma}\step[1.5]\object{X}\\
\hh\step[-.5]\hcu\step[1.5]\id\\
\x\\
\x\\
\hh\id\step\frabox{p}\\
\hh\hcu\step\id\\
\hstep\object{\Gamma}\step[1.5]\object{X}
\end{tangle}
\quad=\quad
\begin{tangle}
\id\step\hx\\
\hx\step\id\\
\hx\step\id\\
\id\step\hx\\
\hh\hcu\step\id\\
\hh\hstep\id\step\frabox{p}\\
\hh\hstep\hcu\step\id
\end{tangle}
\quad=\quad
\begin{tangle}
\hxx\step\id\\
\id\step\hx\\
\id\step\hx\\
\hx\step\id\\
\id\step\hx\\
\id\step\hx\\
\hh\id\step\id\step\frabox{p}\\
\id\step\hddcu\step\id\\
\hh\hcu\Step\id
\end{tangle}
} \\
 && =
\begin{tangle}
\hxx\step\id\\
\id\step\hx\\
\id\step\hx\\
\hx\step\id\\
\id\step\O{\del X_1}\step\id\\
\hh\id\step\id\step\frabox{p}\\
\id\step\hddcu\step\id\\
\hh\hcu\Step\id
\end{tangle}
=\quad
\begin{tangle}
\hxx\step\id\\
\id\step\hx\\
\id\step\hx\\
\O{\del X_1}\step\id\step\id\\
\hh\id\step\id\step\frabox{p}\\
\id\step\hddcu\step\id\\
\hh\hcu\Step\id
\end{tangle}
\quad=\quad\quad
\begin{tangle}
\hxx\step\id\\
\id\step\O{\del X_1}\step\id\\
\O{\del X_1}\step\id\step\id\\
\hh\id\step\id\step\frabox{p}\\
\id\step\hddcu\step\id\\
\hh\hcu\Step\id
\end{tangle}
\quad=\quad
\begin{tangle}
\id\step\O{\del X_1}\step\id\\
\O{\del X_1}\step\id\step\id\\
\hh\id\step\id\step\frabox{p}\\
\id\step\hddcu\step\id\\
\hh\hcu\Step\id
\end{tangle}
\eean
Here the first, fourth and sixth equality are due to (\ref{eq-twist}) and the fifth and
seventh due to the cocommutativity of $\Delta$. Taking the partial trace over $X$ we
obtain 
\[ Tr_{\hat{\2C}}(p)\ \Delta\circ\del X_1= Tr_{\hat{\2C}}(p)\ \del X_1\otimes\del X_1 
   \mcirc\Delta \]
and thus $\Delta\mcirc\del X_1= \del X_1\otimes\del X_1\mcirc\Delta$ since 
$Tr_{\hat{\2C}}p=d(X_1)\ne 0$. Thus $\del X_1\in G=\End(\Gamma,m,\eta)$ is an endomorphism
of the monoid $(\Gamma,m,\eta)$, and by Remark \ref{rem-aut} $G$ is a group.
\qed

\bdefin
An object of $\2C\rtimes\2S$ is homogeneous if there exist $g\in G$ and simple objects
$X_i\in\2C\rtimes\2S, \ i\in\Delta$ such that $X\cong\oplus_i X_i$ and $\del X_i=g$. 
\edefin

\blemma \label{lem-hom}
Let $Z\in\2C\rtimes\2S$ be homogeneous of grade $g$. Then $g$ is still given by
(\ref{e3}). If $(X,p), (Y,q)$ are homogeneous and $(X,p)\cong(Y,q)$ then
$\del(X,p)=\del(Y,q)$.
\elemma

\prf Let $Z\cong\oplus X_1$, where the $X_i$ are simple and $\del X_i=g$. Reviewing the
considerations preceding (\ref{e1}) one sees that this equation remains valid with $X_1$
replaced by $Z$. Thus also (\ref{e3}) holds for homogeneous $Z$, and this is all that is
used in the proof of Proposition \ref{prop_grade}. That isomorphic homogeneous objects
have the same grade is obvious from the definition.
\qed

\bprop \label{prop_cov}
Let $X_1=(X,p)\in\2C\rtimes\2S$ be homogeneous. Then 
$\del\gamma_g(X_1)=g\,\del X_1\,g^{-1}$ for every $g\in G$.
\eprop

\prf Recall that $\gamma_g((X,p))=(X,\gamma_g(p))=(X,p\mcirc g^{-1}\otimes\id_X)$. 
Thus
\[ d(X_1)\,\del\gamma_g(X_1)=\quad
\begin{tangle}
\object{\Gamma}\\
\hh\id\Step\coev\\
\x\step\id\obj{\ol{X}}\\
\x\step\id\\
\hh\id\step\frabox{p}\step\id\\
\id\step\O{g^{-1}}\step\id\step\id\\
\hh\step[-.5]\obj{\Delta}\hstep\hcu\hstep\obj{X}\hstep\hev\\
\hh\hstep\id\\
\hstep\object{\Gamma}
\end{tangle}
\quad=\quad
\begin{tangle}
\object{\Gamma}\\
\hh\id\Step\coev\\
\x\step\id\obj{\ol{X}}\\
\x\step\id\\
\hh\id\step\frabox{p}\step\id\\
\O{g}\step\id\step\id\step\id\\
\hh\step[-.5]\obj{\Delta}\hstep\hcu\hstep\obj{X}\hstep\hev\\
\hstep\O{g^{-1}}\\
\hh\hstep\id\\
\hstep\object{\Gamma}
\end{tangle}
\quad=\quad
\begin{tangle}
\object{\Gamma}\\
\O{g}\\
\hh\id\Step\coev\\
\x\step\id\obj{\ol{X}}\\
\x\step\id\\
\hh\id\step\frabox{p}\step\id\\
\hh\step[-.5]\obj{\Delta}\hstep\hcu\hstep\obj{X}\hstep\hev\\
\hstep\O{g^{-1}}\\
\hh\hstep\id\\
\hstep\object{\Gamma}
\end{tangle}
\quad=d(X_1)\,g\,\del X_1\, g^{-1}.
\]
Here we have used the equation
$\id_\Gamma\otimes g^{-1}\mcirc\Delta=g\otimes\id_\Gamma\mcirc\Delta\mcirc g^{-1}$ which
follows from $\Delta\mcirc g=g\otimes g\mcirc\Delta$.
\qed

The following definition is a variant of a notion due to Turaev \cite{t2}.

\bdefin \label{def_crossed} 
Let $G$ be a (discrete) group. A strict crossed G-category is a strict tensor category
$\2D$ together with 
\begin{itemize}
\item a full tensor subcategory $\2D_G\subset\2D$ of homogeneous objects,
\item a map $\del: \Obj\,\2D_G\rarr G$ constant on isomorphism classes,
\item a (strict) homomorphism $\gamma: G\rarr\Aut\,\2D$. (Here $\Aut\,\2D$ is the group of
invertible strict tensor functors $\2D\rarr\2D$ respecting the braiding.)
\end{itemize}
such that 
\begin{enumerate}
\item $\del(X\otimes Y)=\del X\,\del Y$ for all $X,Y\in\2D_G$. 
\item $\gamma_g(\2D_h)\subset\2D_{ghg^{-1}}$, where $\2D_g\subset\2D_G$ is the full
subcategory $\del^{-1}(g)$.
\end{enumerate}
If $\2D$ is additive we require that every object of $\2D$ be a direct sum of objects in 
$\2D_G$. 
\edefin

\brem
1. A map $\del: \Obj\2D_G\rarr G$ constant on iso-classes and satisfying
$\del(X\otimes Y)=\del(X)\del(Y)$ is the same as a tensor functor $\2D_G\rarr\2G$, where
$\2G$ is the discrete strict monoidal category with $\Obj\2G=G$.

2. In \cite{t2}, $\2D_G=\2D$ was assumed. Since we are working with
additive categories, in particular having all finite direct sums, we must allow
inhomogeneous objects. This added generality will be important later on.

3. Obviously, the definition can be generalized to non-strict tensor categories,
cf.\ \cite{t2}. Also the G-action can be generalized by relaxing the $\gamma_g$ to be 
self-equivalences satisfying natural isomorphisms $\gamma_g\gamma_h\cong\gamma_{gh}$ with
suitable coherence, cf.\ e.g.\ \cite[p.238]{carr}. For our purposes, in particular the
application to conformal field theory \cite{mue15}, the above strict version is sufficient. 
\erem

In view of Definition \ref{def_crossed}, Propositions \ref{prop_grade}, \ref{prop_cov}
essentially amount to the following statement.

\bprop
$\2C\rtimes\2S$ is a crossed G-category, where $\2S\simeq\Rep\,G$.
\eprop

\prf We define $(\2C\rtimes\2S)_G\subset\2C\rtimes\2S$ to be the full subcategory of
homogeneous objects, and we extend $\del$ to $(\2C\rtimes\2S)_G$ in the obvious fashion. 
We have already defined an action $\gamma$ of $G$ on $\2C\rtimes\2S$. Now property 2
follows from Proposition \ref{prop_cov}, but property 1 requires proof. Thus let 
$(X,p),(Y,q)\in\2C\rtimes\2S$ be homogeneous. In view of Lemma \ref{lem-hom} we may
compute
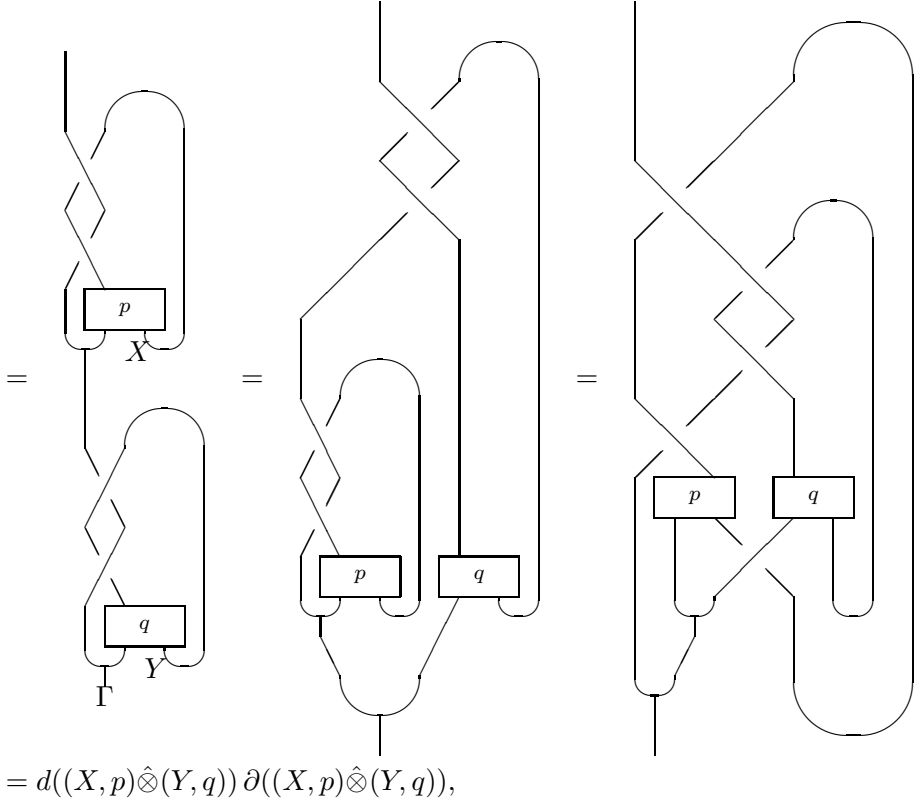
\bean \lefteqn{ d(X,p)d(Y,q)\,\del(X,p)\del(Y,q)= } \\
&& = \quad
\begin{tangle}
\id\step\coev\\
\hx\Step\id\\
\hx\Step\id\\
\hh\id\step\frabox{p}\step\id\\
\hh\hcu\hstep\obj{X}\hstep\hev\\
\hstep\id\step\coev\\
\hstep\hxx\Step\id\\
\hstep\hxx\Step\id\\
\hh\hstep\id\step\frabox{q}\step\id\\
\hh\hstep\hcu\hstep\obj{Y}\hstep\hev\\
\step\object{\Gamma}
\end{tangle}
\quad=\quad
\begin{tangle}
\Step\id\Step\coev\\
\Step\x\Step\id\\
\Step\x\Step\id\\
\step\ne2\Step\id\Step\id\\
\id\step\coev\step\id\Step\id\\
\hx\Step\id\step\id\Step\id\\
\hx\Step\id\step\id\Step\id\\
\hh\id\step\frabox{p}\step\id\step\frabox{q}\step\id\\
\hdcu\hstep\hstep\hev\dd\step\hev\\
\step\cu
\end{tangle}
\quad=\quad
\begin{tangle}
\id\step[4]\mcoev\\
\id\step[3]\ne2\step[3]\id\\
\x\Step\coev\step\id\\
\id\Step\x\Step\id\step\id\\
\id\Step\x\Step\id\step\id\\
\x\Step\id\Step\id\step\id\\
\hh\id\step\frabox{p}\Step\frabox{q}\step\id\step\id\\
\id\step\id\step\xx\step\id\step\id\step\id\\
\id\step\hddcu\Step\id\step\hev\step\id\\
\hcu\step[3]\mev
\end{tangle}
\\
&&= d((X,p)\hat{\otimes}(Y,q))\,\del((X,p)\hat{\otimes}(Y,q)),
\eean
which is the desired result.
\qed

\bdefin \label{def_braid}
A braiding for a crossed G-category $\2D$ is a family of isomorphisms 
$c_{X,Y}: X\otimes Y\rarr{}^X\!Y\otimes X$, defined for all 
$X\in\2D_G,\ Y\in\2D$, such that
\[ \begin{diagram}
  X\otimes Y & \rTo^{s\otimes t} & X'\otimes Y' \\
  \dTo_{c_{X,Y}} && \dTo^{c_{X',Y'}} \\
  {}^X\!Y\otimes X & \rTo_{{}^X\!t\otimes s} & {}^{X'}\!Y'\otimes X'
\end{diagram} \]
commutes for all $s: X\rarr X', t: Y\rarr Y'$, and
\bea c_{X,Z\otimes T} &=& \id_{{}^X\!Z}\otimes c_{X,T}\mcirc c_{X,Z}\otimes\id_T, 
        \label{eq-br1}\\
   c_{X\otimes Y,Z} &=& c_{X,{}^Y\!Z}\otimes\id_Y\mcirc\id_X\otimes c_{Y,Z}, \label{eq-br2}
\eea
for all $X,Y\in\2D_G, \ Z,T\in\2D$
\edefin

\brem 
Motivated by applications to algebraic topology (rather than 3-manifolds as in \cite{t2}),
a special class of braided crossed G-categories was introduced independently in
\cite[Definition 2.1]{carr}. The `categorical G-crossed modules' considered there are
braided crossed G-categories that are also categorical groups, i.e.\ monoidal groupoids
whose objects are invertible up to isomorphism w.r.t.\ $\otimes$. 
\erem

\btheor \label{t-main1}
$\2C\rtimes\2S=\hat{\2C}$ is a braided crossed $G$-category, where $\2S\simeq \Rep\,G$. 
\etheor

\prf Let $X_1=(X,p)\in\hat{\2C}_G,\ Y,Z\in\2C$
and $s\in\Hom_{\hat{\2C}}(Y,Z)=\Hom_\2C(\Gamma\otimes Y,Z)$. We calculate
\bean \iota(c_{X,Z})\ \hat{\circ}\ p\hat{\otimes} s & = & \quad
\begin{tangle}
\step\object{Z}\Step\object{(X,p)}\\
\step\xx\\
\frabox{p}\Step\frabox{s}\\
\step\id\Step\id\\
\step\object{(X,p)}\Step\object{Y}
\end{tangle}
\quad\quad \hat{=}  \quad\quad
\begin{tangle}
\step\object{Z}\Step\object{X}\\
\step\xx\\
\frabox{p}\Step\frabox{s}\\
\id\step\xx\step\id\\
\hh\hcu\Step\id\step\id\\
\hstep\object{\Gamma}\step[2.5]\object{X}\step\object{Y}
\end{tangle}
\quad\quad = \quad\quad
\begin{tangle}
\step[1.5]\object{Z}\step[1.5]\object{X}\\
\hh\step[1.5]\id\step[1.5]\id\\
\hh\step\frabox{s}\step\id\\
\step\id\step\hxx\\
\step\hxx\step\id\\
\hh\frabox{p}\step\id\step\id\\
\id\step\hxx\step\id\\
\hh\hcu\step\id\step\id\\
\hstep\object{\Gamma}\step[1.5]\object{X}\step\object{Y}
\end{tangle} \\
& =& 
\begin{tangle}
\hstep\object{Z}\step[1.5]\object{X}\\
\hh\hstep\id\step[1.5]\id\\
\hh\frabox{s}\step\id\\
\id\step\hxx\\
\hxx\step\d\\
\hxx\Step\id\\
\hh\id\step\frabox{p}\step\id\\
\hh\hcu\step\id\step\id\\
\hstep\object{\Gamma}\step[1.5]\object{X}\step\object{Y}
\end{tangle} 
\quad\quad = \quad\quad
\begin{tangle}
\step[2.5]\object{Z}\step[2.5]\object{X}\\
\hh\step[2.5]\id\step[2.5]\id\\
\hh\Step\frabox{s}\Step\id\\
\step\dd\step\xx\\
\hh\step\frabox{(\del X_1)^{-1}}\step\id\Step\id\\
\hh\step\id\Step\id\Step\id\\
\hh\step\id\step\frabox{p}\Step\id\\
\hh\step\hcu\step\id\Step\id\\
\step[1.5]\object{\Gamma}\step[1.5]\object{X}\Step\object{Y}
\end{tangle} 
\quad = \quad\quad
\begin{tangle}
\hstep\object{Z}\step[2.5]\object{X}\\
\hh\hstep\id\step[2.5]\id\\
\hh\frabox{\gamma_{\del X_1}(s)}\Step\id\\
\id\step\xx\\
\hh\id\step\frabox{p}\step\id\\
\hh\id\step\id\step\id\step\id\\
\hh\hcu\step\id\step\id\\
\hstep\object{\Gamma}\step[1.5]\object{X}\step\object{Y}
\end{tangle} 
\quad\quad\\  \\ \\
&\hat{=}& \quad
\begin{tangle}
\hstep\object{Z}\Step\object{(X,p)}\\
\hh\hstep\id\Step\id\\
\hh\frabox{\gamma_{\del X_1}(s)}\step[1.5]\id\\
\hstep\xx\\
\hh\frabox{p}\step[1.5]\id\\
\hh\hstep\id\Step\id\\
\hstep\object{(X,p)}\Step\object{Y}
\end{tangle}
\quad = \
\gamma_{\del X_1}(s)\hat{\otimes}\id_{X_1} \ \hat{\circ}\ \iota(c_{X,Z})\ \hat{\circ}\
   p\hat{\otimes}\id_Y. 
\eean
We have used cocommutativity of $\Delta$, eq.\ (\ref{eq-twist}), Proposition
\ref{prop_grade} according to which $\del(X,p)\in G$, and the definition of
$\gamma_g\in\Aut\,\2C\rtimes\2S$.  

Let now $(X,p)\in(\2C\rtimes\2S)_g$ and $(Y,q)\in\2C\rtimes\2S$. 
Then the above computation and Lemma \ref{lem-nat1} imply
\[ \iota(c_{X,Y})\ \hat{\circ}\ p\hat{\otimes} q
   =\gamma_g(q)\hat{\otimes} p\ \hat{\circ}\ \iota(c_{X,Y}), \]
thus this expression defines an isomorphism 
$c_{(X,p),(Y,q)}\in\Hom_{\hat{\2C}}((X,p)\otimes(Y,q),\gamma_g(Y,q)\otimes(X,p))$. By
definition, the family $(c_{(X,p),(Y,q)})$ it is natural in the sense of Definition
\ref{def_braid}. The straightforward verification of the braid relations
(\ref{eq-br1}-\ref{eq-br2}) is omitted. 
\qed

\subsection{The $G$-Spectrum of a Galois Extension}

\bdefin \label{def_spec}
The $G$-spectrum $\mathrm{Spec}\,\2D$ of a G-crossed category $\2D$ is set 
$\{g\in G\ |\ \2D_g\ne\emptyset\}$. The $G$-spectrum of a crossed G-category is full if
it coincides with $G$ and trivial if it is $\{e\}$.
\edefin

\blemma
The $G$-spectrum of a crossed G-category $\2D$ contains the unit, is closed under
multiplication and under conjugation with elements of $G$. It is closed under inverses if
$\2D$ has duals, in which case $\mathrm{Spec}\,\2D$ is a normal subgroup of $G$.
\elemma

\prf The first sentence follows from requirements 1 and 2 in Definition \ref{def_crossed}
and the second from the fact that $\del\,\ol{X}=(\del\,X)^{-1}$, which follows from
$\11\prec X\otimes\ol{X}$.
\qed

\bprop
Let $\2D$ be a semisimple rigid crossed $G$-category. Defining $\dim\2D_g$ to be the sum
over the squared dimensions of the simple objects of grade $g$, we have
\[ \dim\2D_g=\dim\2D_e\ \quad\forall g\in\mathrm{Spec}\,\2D. \]
\eprop

\prf Let $\Delta_e,\Delta_g$ be the sets of iso-classes of simple objects in
$\2D_e,\2D_g$, respectively, and let $\{X_i, i\in\Delta_e\}$ and $\{Y_j,j\in\Delta_g\}$ be
representing objects. For $g\in\mathrm{Spec}\,\2D$ we may pick a simple object
$Z\in\2D_g$, and in view of $X_i\otimes Z\in\2D_g$ we have
\bean d(Z)\sum_{i\in\Delta_e} d(X_i)^2 &=& \sum_{i\in\Delta_e} d(X_i)d(X_i\otimes Z)
   = \sum_{i\in\Delta_e}\sum_{j\in\Delta_g} d(X_i)d(Y_j)\dim\Hom(X_i\otimes Z,Y_j) \\
   &=& \sum_{i\in\Delta_e}\sum_{j\in\Delta_g} d(X_i)d(Y_j)\dim\Hom(X_i,Y_j\otimes\ol{Z})
   =\sum_{j\in\Delta_g} d(Y_j)d(Y_j\otimes\ol{Z}) \\
   &=&d(\ol{Z})\sum_{j\in\Delta_g} d(Y_j)^2. \eean
Since $d(Z)=d(\ol{Z})\ne 0$, the claim follows.
\qed

\bprop \label{prop_e}
Let $\2C,\2S$ be as in the preceding section. The embedding
$(\2C\cap\2S')\rtimes\2S\hookrightarrow\2C\rtimes\2S$ gives rise to an isomorphism
$(\2C\rtimes\2S)_e\cong(\2C\cap\2S')\rtimes\2S$. $\2C\rtimes\2S$ has trivial $G$-spectrum iff
$\2S\subset Z_2(\2C)$.
\eprop

\prf If $X\in\2C\cap\2S'$ then $c_{X,\Gamma}\circ c_{\Gamma,X}=\id$, thus every simple
summand of $\iota(X)$ has grade $e$. This implies 
$(\2C\cap\2S')\rtimes\2S\subset(\2C\rtimes\2S)_e$. As to the converse, every simple object
$X_1\in\2C\rtimes\2S$ is isomorphic to one of the form $(X,p)$, where $X\in\2C$ is
simple and $p$ is a minimal idempotent. In \cite[Proposition 4.2]{mue06} it was shown that
the action $\gamma$ of $G$ on $\2C\rtimes\2S$ acts transitively on the minimal central
idempotents in $\End_{\hat{\2C}}(\iota(X))$, in particular all simple summands of
$\iota(X)$ appear with the same multiplicity $N$. If $\iota(X)\cong N\oplus_i (X,p_i)$ is
the decomposition into simples, we conclude from Proposition \ref{prop_cov} that the set
$\{ \del(X,p_i)\}$ is a conjugacy class in $G$. If $X_1\prec\iota(X)$ has grade $e$ then
this conjugacy class is $\{e\}$, thus $\del(X,p_i)=e$ for all $i$. This means
\begin{equation} \label{e-xx}
\begin{tangle}
\object{\Gamma}\\
\hh\id\Step\coev\\
\x\step\id\obj{\ol{X}}\\
\x\step\id\\
\hh\id\step\frabox{p_i}\step\id\\
\hh\step[-.5]\obj{\Delta}\hstep\hcu\hstep\obj{X}\hstep\hev\\
\hh\hstep\id\\
\hstep\object{\Gamma}
\end{tangle}
\quad=\quad
\begin{tangle}
\hh\id\step\step\coev\\
\hh\id\step\frabox{p_i}\step\id\\
\id\step\obj{\eta}\counit\step\hev\\
\object{\Gamma}
\end{tangle}
\end{equation}
for all minimal central idempotents $p_i$ in $\End_{\hat{\2C}}\iota(X)$. By linearity,
(\ref{e-xx}) holds for all central idempotents, in particular for
$\id_{\iota(X)}=\ve\otimes\id_X$. Plugging this into (\ref{e-xx}) we obtain
$(\id\otimes Tr_X)(c_{X,\Gamma}\circ c_{\Gamma,X})=d(X)\id_\Gamma$, and by naturality we
conclude
\[ S(X,Y)=(Tr_Y\otimes Tr_X)(c_{X,Y}\circ c_{Y,X})=d(X)d(Y) \]
for all simple $Y\in\2S$. By \cite[Proposition 2.5]{mue11} this is equivalent to
$X\in\2C\cap\2S'$. Now, triviality of the $G$-spectrum is equivalent to 
$\2C\rtimes\2S=(\2C\rtimes\2S)_e=(\2C\cap\2S')\rtimes\2S$, which in turn is equivalent
to $\2C\cap\2S'=\2C$ and finally to $\2S\subset Z_2(\2C)$.
\qed

\brem We emphasize one observation made in the proof: Whereas every simple object $X_1$ of
$\2C\rtimes\2S$ defines an element $\del X_1$ of $G$, every simple object $X\in\2C$
defines a unique conjugacy class in $G$.
\erem

Let $\2S_0\subset\2S$ be a full subcategory, where both categories are finite admissible
STCs. Let $(\Gamma,\ldots)$, $(\Gamma_0,\ldots)$ be the corresponding Frobenius algebras
in $\2S_0, \2S$, respectively, with automorphism groups $G_0, G$. Then
$\Gamma\cong\Gamma_0\oplus Z$ and $\Hom(\Gamma_0,Z)=\{0\}$, thus the projector
$q\in\End\,\Gamma$ onto $\Gamma_0$ is central. The group 
\[ N=\{ g\in G\ | \ g\circ q=q\} \]
is a normal subgroup of $G=\Aut(\Gamma,m,\eta)$. It coincides with
\[ N=\{g\in G\ | \ \pi_X(g)=\id_{E(X)}\ \forall X\in\2S_0 \}, \]
where $E: \2S\rarr\mathrm{Vect}_\7C$ is the fiber functor and $\pi_X$ is the
representation of $G$ on $E(X)$. This is easily deduced from 
$E(X)=\Hom(\11,\Gamma\otimes X)$ and the fact that $g\in G$ acts on $E(X)$ by
$\pi_X(g): \phi\mapsto g\otimes\id_X\mcirc\phi$. This implies $G_0\cong G/N$.

\btheor \label{theor-spec}
Let $\2S\subset\2C$ with $\2S\simeq\Rep\,G$. Then $\mathrm{Spec}\,\2C\rtimes\2S=N$, where
$N$ is the normal subgroup of $G$ corresponding to the full inclusion 
$\2S\cap Z_2(\2C)\subset\2S$ as above. $\2C\rtimes\2S$ has full $G$-spectrum iff 
$\2S\cap Z_2(\2C)$ is trivial, i.e.\ consists only of multiples of $\11$.
\etheor

\prf Let $q\in\End_{\hat{\2C}}(\Gamma)$ be the projection onto $\Gamma_0$, and let 
$v: \Gamma_0\rarr\Gamma, v': \Gamma\rarr\Gamma_0$ satisfy $v\circ v'=q$, 
$v'\circ v=\id_{\Gamma_0}$. Then with $X_1=(X,p)\in(\2C\rtimes\2S)_G$ we have
\[ d(X_1)\ q\mcirc\del(X,p)=\quad
\begin{tangle}
\O{q}\\
\hh\id\Step\coev\\
\x\step\id\obj{\ol{X}}\\
\x\step\id\\
\hh\id\step\frabox{p}\step\id\\
\hh\step[-.5]\obj{\Delta}\hstep\hcu\hstep\obj{X}\hstep\hev\\
\hh\hstep\id\\
\hstep\object{\Gamma}
\end{tangle}
\quad=\quad
\begin{tangle}
\O{v}\\
\hh\id\Step\coev\\
\x\step\id\obj{\ol{X}}\\
\step[-.5]\mobj{\Gamma_0}\hstep\x\step\id\\
\O{v'}\Step\id\step\id\\
\hh\id\step\frabox{p}\step\id\\
\hh\step[-.5]\obj{\Delta}\hstep\hcu\hstep\obj{X}\hstep\hev\\
\hh\hstep\id\\
\hstep\object{\Gamma}
\end{tangle}
\quad=\quad
\begin{tangle}
\O{v}\\
\hh\id\Step\coev\\
\O{v'}\Step\id\step\id\\
\hh\id\step\frabox{p}\step\id\\
\hh\step[-.5]\obj{\Delta}\hstep\hcu\hstep\obj{X}\hstep\hev\\
\hh\hstep\id\\
\hstep\object{\Gamma}
\end{tangle}
\quad
=d(X_1)q,
\]
where we used $\Gamma_0\in\2S\cap Z(\2C)$. We conclude
$\mathrm{Spec}\,\2C\rtimes\2S\subset N$. 

In a braided crossed G-category $\2D$ we have isomorphisms 
$c_{X,Y}: X\otimes Y\rarr\gamma_g(Y)\otimes X$ whenever $X\in\2D_G$. By definition, 
$g\in\mathrm{Spec}\,\2D$, thus in the fixpoint category $\2D^{\mathrm{Spec}\,\2D}$ the
action $\gamma_g$ disappears and $\2D^{\mathrm{Spec}\,\2D}$ is braided in the usual
sense. We therefore have an intermediate extension
\[ \2C\ \subset\ (\2C\rtimes\2S)^{\mathrm{Spec}\,\2C\rtimes\2S} \ \subset\ \2C\rtimes\2S
\]
that is braided. On the other hand, in view of Proposition \ref{prop_e} it is clear that
the maximal intermediate extension of $\2C$ that is braided is given by
\[ \2C\ \subset\ \2C\rtimes(\2S\cap Z_2(\2C)) \ \subset\ \2C\rtimes\2S. \]
By the Galois correspondence established in \cite[Section 4.2]{mue06} we have
$\2C\rtimes(\2S\cap Z_2(\2C))=(\2C\rtimes\2S)^N$, where $N$ is as defined above. Now the
inclusion 
\[ (\2C\rtimes\2S)^{\mathrm{Spec}\,\2C\rtimes\2S}\subset 
  \2C\rtimes(\2S\cap Z_2(\2C))=(\2C\rtimes\2S)^N \]
implies $N\subset \mathrm{Spec}\,\2C\rtimes\2S$. This completes the proof of
$\mathrm{Spec}\,\2C\rtimes\2S=N$. The last claim is immediate.
\qed

The following corollary will be very useful in conformal field theory \cite{mue15}.

\bcoro \label{coro-spec}
If $\2C$ is modular and $\Rep\,G\simeq\2S\subset\2C$ then $\2C\rtimes\2S$ has full
$G$-spectrum and $(\2C\rtimes\2S)_e$ is modular.
\ecoro

\prf Modularity of $\2C$ is equivalent to triviality of $Z_2(\2C)$, thus the last
statement of Theorem \ref{theor-spec} implies $\mathrm{Spec}\,\2C\rtimes\2S=G$. 
Since $\2C$ is modular, \cite[Corollary 3.6]{mue11} implies $Z_2(\2C\cap\2S')=\2S$. Thus 
$(\2C\cap\2S')\rtimes\2S$ is modular by \cite[Theorem 4.4]{mue06} and coincides with
$(\2C\rtimes\2S)_e$ by Proposition \ref{prop_e}. 
\qed


\subsection{Abelian Case} \label{ss-34}
Let $X\in\2C$ be simple and let $X_j\in\2C\rtimes\2S,\ j\in J,$ be simple objects such that
$\iota(X)\cong\oplus_{j\in J} X_j$. In \cite{mue06} it was shown that $G$ acts ergodically
on the center of the algebra $\End\,\iota(X)$. In view of $\del\gamma_g(X)=g\del(X)g^{-1}$
this clearly implies that the set $\{ \del X_j\ | \ j\in J\}$ is a conjugacy class in $G$. We
thus obtain a map $\del_0$ from the simple objects in $\2C$ to the conjugacy classes in
$G$.
In the case where $G$ is abelian, all simple summands of $\iota(X)$ have the same grade,
which induces a $G$-grading on the category $\2C$. In the remainder of this subsection we
will give a more explicit description of this grading. 

Let thus $G$ be abelian and $K=\hat{G}$. Then $\Gamma\cong\oplus_{k\in K}X_k$, where all
$X_k, k\in K$ are invertible, and
$\End\,\Gamma\cong\oplus_{k\in K}\End\,X_k\cong\oplus_{k\in K}\7F$. By our normalization
$\ve\circ\eta=1$, $p_e=\eta\circ\ve\in\End\,\Gamma$ is an idempotent, projecting on the
summand $X_e$. Let $X\in\2C$ and $(X,p)\in\2C\rtimes\2S$ be simple. By the above
considerations, $\iota(X)$ is homogeneous, thus (\ref{e3}) defines an element of
$\Aut(\Gamma,m,\eta)\cong G$. In view of $X_k\otimes X_l\cong X_{kl}$ we may insert $p_e$
into (\ref{e3}) at the appropriate place, obtaining
\[
\del(X,p)=\quad
\left( \ \ \
\begin{tangle}
\hh\step\coev\\
\hh\frabox{p}\step\id\\
\obj{\eta}\counit\step\hev\\
\end{tangle}
\ \ \right)^{-1}
\cdot\quad
\begin{tangle}
\object{\Gamma}\\
\hh\id\Step\coev\\
\x\step\id\obj{\ol{X}}\\
\x\step\id\\
\hh\id\step\frabox{p}\step\id\\
\id\step\obj{\eta}\counit  \step\hev\\
\id\step[1.2]\obj{\ve}\step[-.2]\unit\\
\hh\step[-.5]\obj{\Delta}\hstep\hcu\\
\hh\hstep\id\\
\hstep\object{\Gamma}
\end{tangle}
\]
Now,
\[ 
\begin{tangle}\hh\step\id\\
\hh\frabox{p}\\
\obj{\eta}\counit\step\id\\
\end{tangle}
\quad=d(X)^{-1}\quad
\begin{tangle}
\hh\step\coev\step\id\\
\hh\frabox{p}\step\id\step\id\\
\obj{\eta}\counit\step\hev\step\id\\
\end{tangle}
\]
and we obtain
\begin{equation} \label{e5}
 \del_0(X)=\del((X,p))=d(X)^{-1}\quad
\begin{tangle}
\hh\id\Step\coev\\
\xx\step\id\\
\xx\obj{X}\step\id\\
\hh\id\Step\ev\\
\object{\Gamma}
\end{tangle} 
\end{equation}

We have thus shown:

\bprop Consider $\2S\subset\2C$ where $\2S$ is symmetric, even and all its simple objects
are one dimensional, equivalently $\2S\simeq\Rep\,G$ with $G$ abelian. Let
$(\Gamma,m,\eta)$ be the regular monoid in $\2S$. Then (\ref{e5}) defines an element
$\del_0X$ of $G$ for every simple $X\in\2C$. If we define $\2C_G$ to be the full
subcategory of homogeneous objects, i.e.\ of objects all simple summands $X_j$ of which
have the same $\del_0 X_j$, then $\2C$ is a $G$-graded tensor category.
(To wit, $\2C$ is a crossed $G$-category in the sense of Definition \ref{def-crossedG}
with trivial $G$-action.)
\eprop

\brem This result can be obtained in a more direct way. It suffices to notice that 
the map $\varphi_X: K\rarr\7F$ defined by 
$\varphi_X(k)\id_{X_k}=(\id_{X_k}\otimes Tr_X)(c_{X,X_k}\circ c_{X_k,X})$ is a character
of $K$, thus an element of $G$. (This goes back at least to \cite{khr1}.) From the above
considerations it is then clear that the two definitions yield the same element
$\del_0X\in G$.
\erem


\vspace{1cm}\noindent
{\it Acknowledgments.} Financial support by NWO is gratefully acknowledged.

\end{document}